\documentclass[12pt,a4paper]{article}

\usepackage{amssymb}
\usepackage{amsfonts}
\usepackage{amsmath}
\usepackage{latexsym,epsfig}
\usepackage{graphicx} 
\usepackage[usenames,dvipsnames]{color} 
\usepackage{breakurl} 
\usepackage{soul} 
\usepackage{array} 
\usepackage{color} 
\usepackage{colortbl} 
\usepackage{verbatim}
\usepackage[hang,small,bf]{caption} 
\usepackage[ruled,linesnumbered]{algorithm2e}  
\usepackage{adjustbox}
\usepackage{url}
\usepackage{enumerate} 
\usepackage{rotating}
\usepackage{hyperref} 

\newcommand{\PSet}{\mathcal{P}}

\newcommand{\Z}{\mathbb{Z}}

\title{Searching for Large Circulant Graphs}

\author{Ramiro Feria-Puron$^{1}$, Hebert P\'erez-Ros\'es$^{2, 1}$, Joe Ryan$^{1}$ \\ \\
$\phantom{0}^{1}$\emph{Dept. of Software Engineering and Computer Science} \\ 
\emph{University of Newcastle, Australia} \\\\
\vspace{2mm} 
$\phantom{0}^{2}$\emph{Dept. of Mathematics, University of Lleida, Spain} } 

\begin{document}
\maketitle

\vskip 5mm
\thispagestyle{empty}

\begin{abstract}
We address the problem of constructing large undirected circulant networks with given degree and diameter. First we discuss the theoretical upper bounds and their asymptotics, and then we describe and implement a computer-based method to find large circulant graphs with given parameters. For several combinations of degree and diameter, our algorithm produces the largest known circulant graphs. We summarize our findings in a table, up to degree 15 and diameter 10, and we perform a statistical analysis of this table, which can be useful for evaluating the performance of our methods, as well as other constructions in the future. 
\end{abstract}

\vskip 5mm

\section{Introduction}
\label{sec:intro}

The Degree$/$Diameter Problem (or DDP, for short) is one of the main theoretical issues in network design. DDP asks for constructing the largest possible graph, in terms of the number of vertices, subject to an upper bound on the maximum degree $\Delta$ and an upper bound on the diameter $D$. Let $N_{\Delta,D}$ be the number of vertices of the largest graph with maximum degree $\Delta$ and diameter $D$; it is well known that  

\begin{equation}
\label{eq:moorebound}
\begin{split}
N_{\Delta,D} \leq M_{\Delta,D} &= 1 + \Delta + \Delta(\Delta -1) + \dots + \Delta(\Delta -1)^{D-1} \\
&= \left\{ 
\begin{array}{ll} 1 + \Delta\ \frac{(\Delta -1)^{D} - 1}{\Delta -2} & \mbox{if $\Delta > 2$} \\ 
2D+1 & \mbox{if $\Delta = 2$} \\
\end{array}\right.
\end{split}
\end{equation}

The number $M_{\Delta,D}$ is called the \emph{Moore bound}, and a graph of order $M_{\Delta,D}$ is called a \emph{Moore graph} \cite{Moore-survey}. Moore graphs are rather scarce: they exist only for a few combinations of $\Delta$ and $D$. For $D=1$ and $\Delta \ge 1$, they are the complete graphs on $\Delta + 1$ vertices. For $D \ge 2$ and $\Delta = 2$, they are the cycles on $2D+1$ vertices. For $D = 2$, Moore graphs exist for $\Delta = 2,3,7$, and possibly $57$. 

The Degree$/$Diameter Problem can also be formulated for digraphs, but in this paper we will only be concerned with the undirected version. On the other hand, we are interested in some particular versions of the problem, namely when the graphs are restricted to a certain class, such as the class of bipartite graphs \cite{Del85,bip09a,bip09b,bip12,bip13}, planar graphs \cite{Gob93,Hell93, Fell95,Fell98,Ti08,Ti12}, vertex-transitive graphs \cite{Mac10,Sia07}, Cayley graphs \cite{Bran98,Haf95,Mac10,Sia07,Ve13}, Cayley graphs of abelian groups \cite{Dou04}, or circulant graphs \cite{Wong74,Muga94,Mar03}. In particular, in this paper we are concerned with circulant graphs, i.e. Cayley graphs of finite cyclic groups. 

For most of the aforementioned graph classes there exist Moore-like upper bounds, which are usually smaller than the Moore bound for general graphs. Most of the research, both in the general version, as well as in these restricted versions of DDP, falls into one of two main categories: 

\begin{enumerate}
\item Proving the non-existence of graphs with order close or equal to the upper bound, or
\item Giving constructions of large graphs, whose order approach the upper bounds as much as possible. 
\end{enumerate}

Research in the second category has been substantially boosted by the compliation of record graph tables, containing the largest known graphs for several combinations of degree and diameter. These tables provide benchmarks to test construction methods and computer search algorithms. Thus, after the first compilation started by Comellas in 1995 \cite{tablecomellas}, other record graph tables have been collected for general graphs, bipartite graphs, planar graphs, and Cayley graphs \cite{ourtables}. 

The study of circulant graphs began in 1970 with Elspas and Turner \cite{Els70}. Coincidentally, it was also Elspas who had formulated the Degree$/$Diameter Problem back in 1964 \cite{Els64}, but apparently he missed the connection between both topics. 

Even though the abelian property of the underlying group prevents abelian Cayley graphs in general (and circulant graphs in particular) to grow as large as their non-abelian counterparts, these graphs have been widely used as topologies for computer networks and parallel computers, due to their other nice properties. Paraphrasing \cite{Dou04}: \lq\lq\ldots the extra structure provided by the groups may provide compensating advantages \ldots, such as good routing algorithms, easy constructibility, and the ability to map common problems onto the architecture\rq\rq. 

However, up to now there were no comprehensive tables for abelian Cayley graphs, or their subfamily of circulant graphs. This is sufficient motivation to start a collection of benchmarks in the class of circulant graphs, which can be used for comparison purposes in the future. Our goals in this paper are twofold: 

\begin{enumerate}
\item Design and implement a computer search algorithm to find large circulant graphs with small degree and diameter, and
\item Start the compilation of a table containing the largest known circulant graphs for some small values of degree and diameter, including the results obtained with the aforementioned algorithm. 
\end{enumerate}

Additionally, we explore the sharpness of the existing upper bounds for some particular combinations of degree and diameter.   

\section{Definitions and basic facts}
\label{sec:def}

An undirected circulant graph $C(n; S)$ is a Cayley graph on the cyclic group $\Z_n$, with a symmetric connection set $S$ (i.e. $S = S^{-1}$). Since $\Z_n$ is abelian, we can switch to additive notation and rephrase the symmetry condition of the connection set as $S = -S$. In order to simplify the notation, we drop the curly brackets $\{ \}$ in the specification of the set $S$ in $C(n; S)$. 

Being a Cayley graph, $C(n; S)$ is vertex-transitive. The degree of $C(n; S)$ is $\Delta = \vert S \vert$, and its order is obviously $n$. A circulant graph can also be defined as a graph of $n$ vertices whose adjacency matrix is circulant \cite{Els64}.

Regarding the degree, we distinguish two cases: 

\begin{enumerate}
\item Even degree: $\Delta = 2t$. In that case, $S = \{ \pm s_1, \ldots, \pm s_t \}$, where $1 \leq s_1 < \ldots < s_t < \frac{n}{2}$.
\item Odd degree: $\Delta = 2t+1$. In that case, $S = \{ \pm s_1, \ldots, \pm s_t, \frac{n}{2} \}$, where $1 \leq s_1 < \ldots < s_t < \frac{n}{2}$. It follows that odd degree is only possible when $n$ is even. 
\end{enumerate}

$C(n; S)$ is connected if, and only if, $\gcd(n, s_1, \ldots, s_t) = 1$. If $\gcd(n, r) = 1$, then $C(n; S)$ is isomorphic to $C(n; rS)$, where multiplication is taken modulo $n$. In that case we say that the connection sets $S$ and $rS$ are multiplicatively related. It should be noted, however, that two circulant graphs may be isomorphic without their connection sets being multiplicatively related \cite{Els64, Muzy04}. 

Now let $N^{circ}_{\Delta,D}$ be the number of vertices of the largest circulant graph with degree $\Delta$ and diameter $D$. It was proved in \cite{BW85} that, if $\Delta = 2t$, then

\begin{equation}
\label{eq:circumoore}
N^{circ}_{\Delta,D} \leq F(t,D) = \sum_{i=0}^{t} 2^i {t \choose i} {D \choose i}
\end{equation}

This upper bound was later rediscovered by Muga \cite{Muga94}. The quantity $F(t,D)$ also turns out to be an upper bound for $N^{AC}_{\Delta,D}$, the order of the largest Cayley graph over an abelian group, with degree $\Delta$ and diameter $D$ \cite{Dou04}. It is quite surprising that no better general upper bound (yet) exists for circulant graphs, considering that they are a special case of abelian Cayley graphs.  

The numbers $F(t,D)$ of Eq. \ref{eq:circumoore} are known as \emph{Delannoy numbers} (sequence A008288 of \cite{oeis}), and they arise in a variety of combinatorial and geometric problems \cite{sulanke}. For example, they correspond to the volume of the ball of radius $D/2$ in the $L^1$ metric in $t$ dimensions \cite{Dou04, mesh12, atanassov}. 

Unaware of the Delannoy connection, Stanton and Cowan had already studied these numbers back in 1970 \cite{stan70}, as a generalization of the binomial coefficients, which satisfy the recurrence:

\begin{equation}
\label{eq:recurrence1}
\begin{split}
F(t,D) &= F(t-1,D) + F(t,D-1) + F(t-1,D-1), \mbox{ with} \\
F(t,1) &= 2t+1, \; \mbox{ for } t \geq 0.
\end{split}
\end{equation}

They gave several interesting formulas for these numbers, such as:  

\begin{equation}
\label{eq:delannoy}
F(t,D) = \sum_{i=0}^{t} {t \choose i} {{D+i} \choose t} = \sum_{i=0}^{t} {{D+i} \choose i} {D \choose {t-i}}
\end{equation}


\begin{figure}[htbp]
\begin{center}
	 	\includegraphics[width=0.95\textwidth]{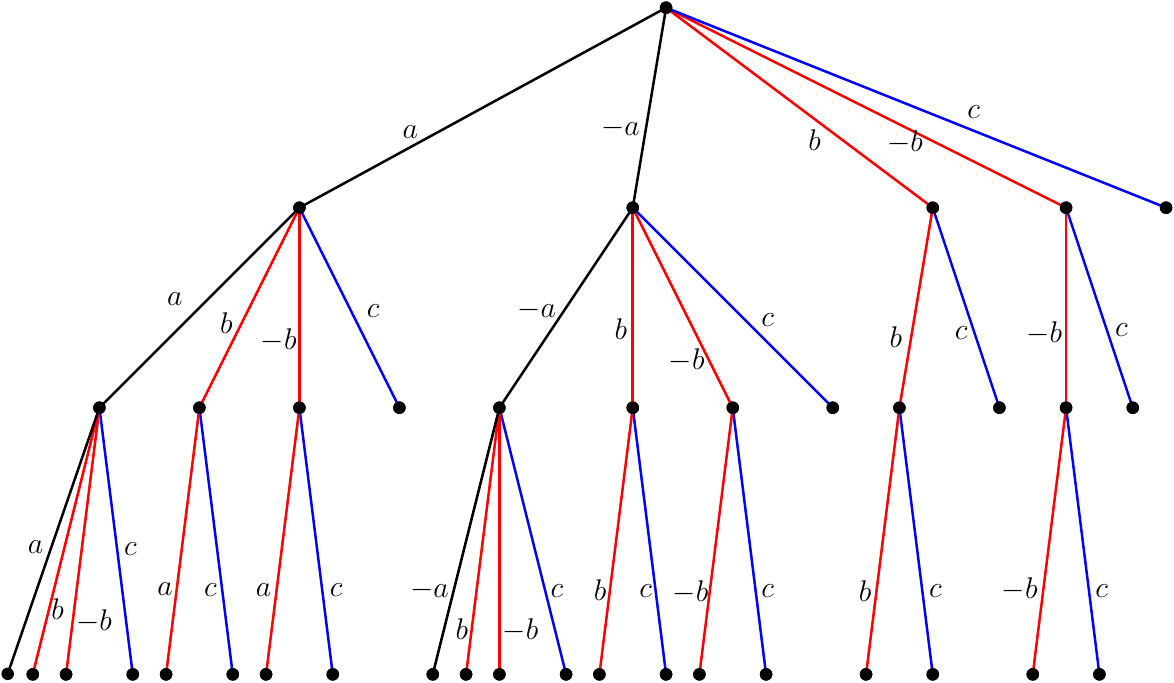}
	 \caption{Tree representation of a maximal abelian Cayley graph of degree $5$ and diameter $3$.}
	 \label{fig:tree}
\end{center}	
\end{figure}

In the case of odd $\Delta$ (i.e. $\Delta = 2t+1$) we have the generator $\frac{n}{2}$, which is its own inverse. Figure \ref{fig:tree} provides a graphical example of that case, for $\Delta = 5$ and $D = 3$. We have denoted the generators as $a, b, c$, where $c = -c$. We can see in this example that the edges labeled with $c$ duplicate every vertex, except those in the lowest level (level $D$). Therefore, an upper bound for the maximum number of vertices in this case is: 

\begin{equation}
\label{eq:oddcase}
F'(t,D) = F(t,D) + F(t,D-1)
\end{equation}

From Equations \ref{eq:recurrence1} and \ref{eq:oddcase} we can deduce that the numbers $F'(t,D)$ also satisfy the same recurrence: 

\begin{equation}
\label{eq:recurrence2}
\begin{split}
F'(t,D) &= F'(t-1,D) + F'(t,D-1) + F'(t-1,D-1), \mbox{ with} \\
F'(t,1) &= 2t+2, \; \mbox{ for } t \geq 0.
\end{split}
\end{equation}

Additionally, $F(t,D) < F'(t,D) < F(t+1,D)$, i.e. the sequence of upper bounds obtained by fixing the diameter is monotonically increasing. Nevertheless, $F(t,D)$ and $F'(t,D)$ grow at slightly different rates, as we will see next. 

In order to find an explicit formula for $F'(t,D)$, we proceed as in \cite{mesh12}: Define the generating function $A_D(z) = \sum_{t \geq 0}{F'(t,D)z^t}$. Multiplying both sides of Eq. \ref{eq:recurrence2} by $z^t$ and summing over $t \geq 1$ we get

\begin{equation}
\label{eq:gen1}
A_D(z) - A_D(0) = zA_D(z) + (A_{D-1}(z)-A_{D-1}(0)) + zA_{D-1}(z)
\end{equation}

whence

\begin{equation}
\label{eq:gen2}
A_D(z) = \frac{1+z}{1-z} A_{D-1}(z)
\end{equation}

With the aid of the boundary condition in Eq. \ref{eq:recurrence2} we get

\begin{equation}
\label{eq:gen3}
A_D(z) = \frac{2(1+z)^{D-1}}{(1-z)^{D+1}}.
\end{equation}

Now, $A_D(z)$ is the product of $(1+z)^{D-1} = \sum_t{D-1 \choose t}z^t$ and $1/(1-z)^{D+1} = \sum_t{D+t \choose t}z^t = \sum_t{D+t \choose D}z^t$. Then, the series of $A_D(z)$ can be obtained as the convolution of the respective factor series, which gives us:  

\begin{equation}
\label{eq:oddbound}
\begin{split}
F'(t,D) &= 2\sum_{i=0}^{t}{{D-1 \choose i}{D+t-i \choose t-i}} \\
        &= 2\sum_{i=0}^{t}{{D-1 \choose i}{D+t-i \choose D-i}} \\
        &= 2\sum_{i=0}^{t}{{D-1 \choose t-i}{D+i \choose i}} 
\end{split}
\end{equation}

Note that Eq. \ref{eq:oddbound} is the same as Eq. 2 of \cite{mesh12}, with the roles of $\Delta$ and $D$ swapped.\footnote{More precisely, $t$ takes the value $p$ (which corresponds to $\lfloor D/2 \rfloor$ in \cite{mesh12}), and $D$ takes the value $k$ (which corresponds to $\Delta/2$).} This symmetry arises as a consequence of the symmetry of Eq. \ref{eq:recurrence2}. In the same manner we can derive other formulas for $F(t,D)$ from Eq. 2 of \cite{mesh12}: 

\begin{equation}
\label{eq:evenbound}
\begin{split}
F(t,D) &= \sum_{i=0}^{t}{{D \choose i}{D+t-i \choose t-i}} \\
        &= \sum_{i=0}^{t}{{D \choose t-i}{D+i \choose i}} 
\end{split}
\end{equation}

With the aid of the methods developed in \cite{zeilberger} it can be shown that the numbers $F(t,D)$ and $F'(t,D)$ do not have a closed form, meaning that they cannot be represented as a linear combination of a fixed number of hypergeometric terms. However, we can obtain asymptotic formulas for them. From Corollary 2.1 of \cite{mesh12}, and using the symmetry just discussed, we get that both $F(t,D)$ and $F'(t,D)$ grow asymptotically as $\displaystyle \frac{(2t)^D}{D!} + o(t^D)$. 

The first few values of $F(t,D)$ and $F'(t,D)$ are collected in Table \ref{tab:circulant1}, up to degree 15 and diameter 10. We can fit an exponential function $e^{f(\Delta, D)}$ to these values, where 

\begin{equation}
\label{eq:fit1}
\begin{split}
f(\Delta, D) = &-0.8015 + 0.5612 \Delta + 0.5433 D -0.05421 \Delta^2 + 0.1499 \Delta D \\
               &-0.09784 D^2 + 0.001574 \Delta^3 -0.001924 \Delta^2 D \\
               &-0.003822 \Delta D^2 + 0.004591 D^3, 
\end{split}
\end{equation}

with $R^2 = 0.9986$. 

\begin{figure}[htbp]
\begin{center}
	 	\includegraphics[width=0.9\textwidth]{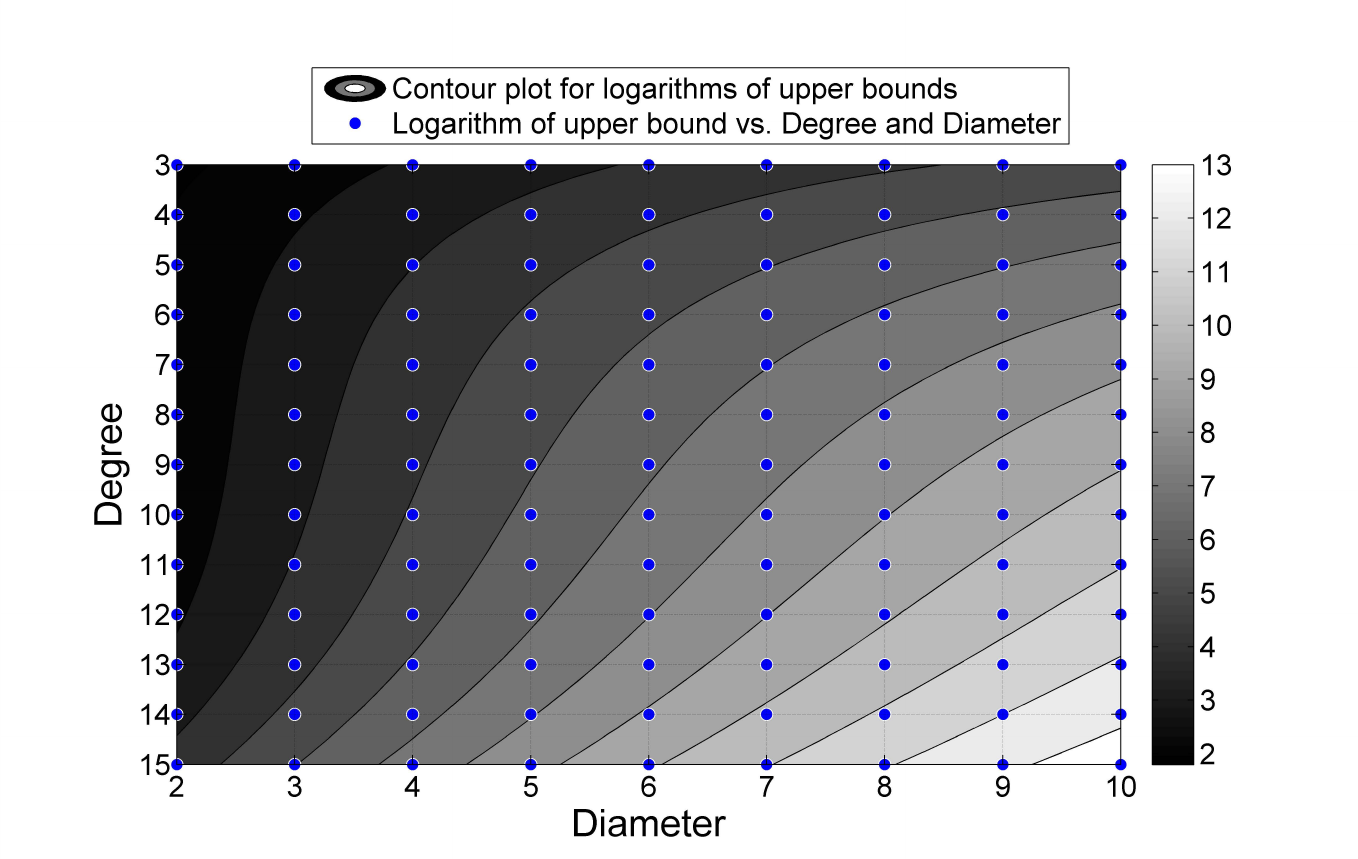}
\end{center}
\caption{Contour plot of the function $f(\Delta, D)$, that approximates the logarithms of the upper bounds $F(t,D)$ and $F'(t,D)$.}
\label{fig:ContourPlotUpperBounds1A}
\end{figure}

Figure \ref{fig:ContourPlotUpperBounds1A} shows a contour plot of $f(\Delta, D)$.

\begin{figure}[htbp]
\begin{center}
	 	\includegraphics[width=0.9\textwidth]{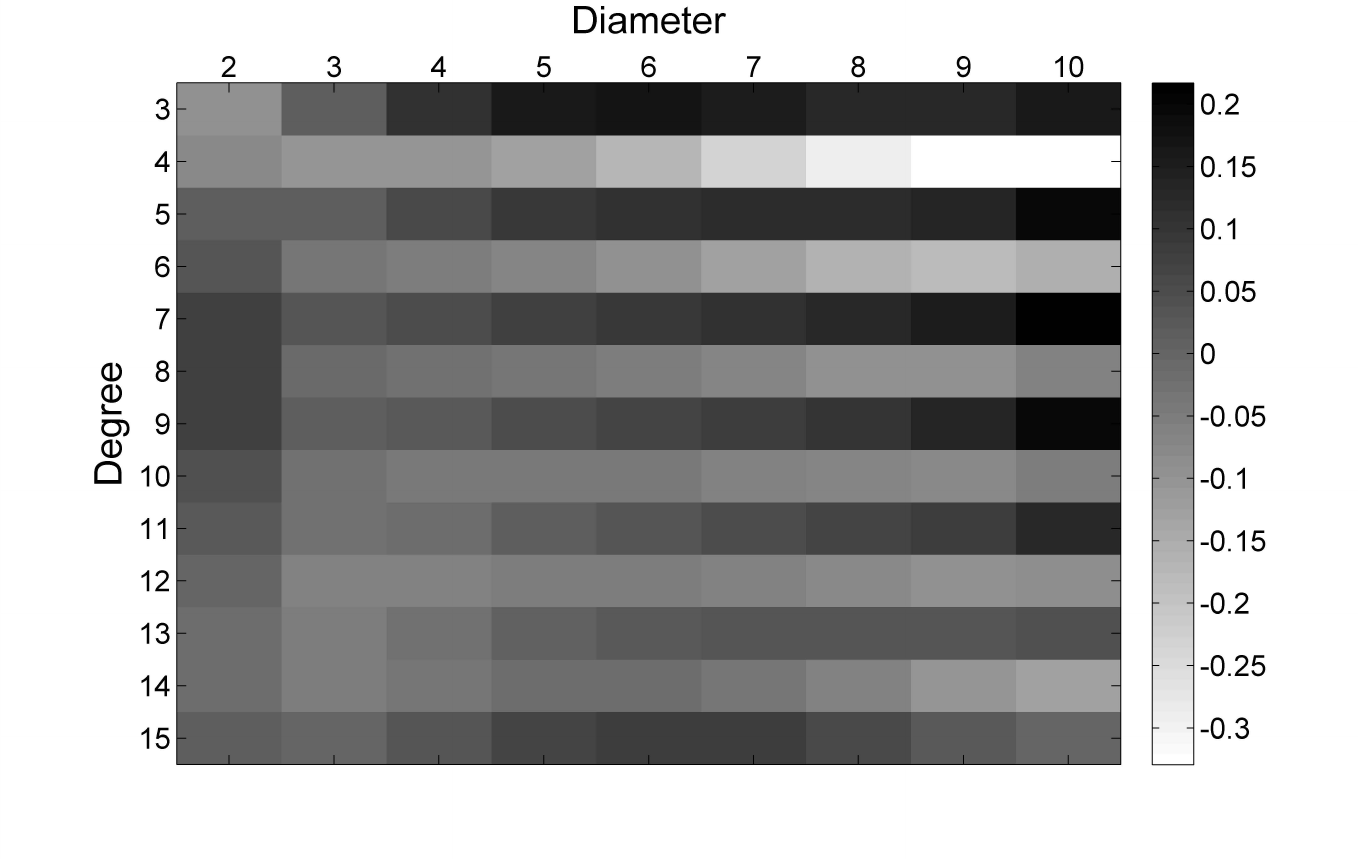}
\end{center}
\caption{Predicted (logarithmic) upper bound minus actual (logarithmic) upper bound. Light entries correspond to negative differences. Darker colors correspond to positive differences.}
\label{fig:Log-upper-BW-1A}
\end{figure}

If we compute the differences between the values predicted by $f(\Delta, D)$ and the logarithms of the actual upper bounds $F(t,D)$ and $F'(t,D)$, we see that the numbers $F(t,D)$ are slightly above $f(\Delta, D)$, while $F'(t,D)$ are slightly below $f(\Delta, D)$, as shown in Figure \ref{fig:Log-upper-BW-1A}. Therefore, it might be useful to compute different approximations for the values with even and odd degree, respectively:  

\begin{equation}
\label{eq:evenfit}
\begin{split}
f(\Delta, D) = &-0.5686 + 0.4266 \Delta + 0.6491 D -0.03746 \Delta^2 \\
               &+ 0.1366 \Delta D -0.1009 D^2 + 0.00098 \Delta^3 \\
               &-0.001492 \Delta^2 D -0.003578 \Delta D^2 + 0.004598 D^3
\end{split}
\end{equation}

(with $R^2 = 0.9999$), and 

\begin{equation}
\label{eq:oddfit}
\begin{split}
f(\Delta, D) = &-0.6153 + 0.5063 \Delta + 0.4993 D -0.04712 \Delta^2 \\
               &+ 0.1507 \Delta D -0.09576 D^2 + 0.001282 \Delta^3 \\
               &-0.001776 \Delta^2 D -0.003974 \Delta D^2 + 0.004584 D^3,
\end{split}
\end{equation}

with $R^2 = 0.997$.

\begin{figure}[htbp]
\begin{center}
	 	\includegraphics[width=0.9\textwidth]{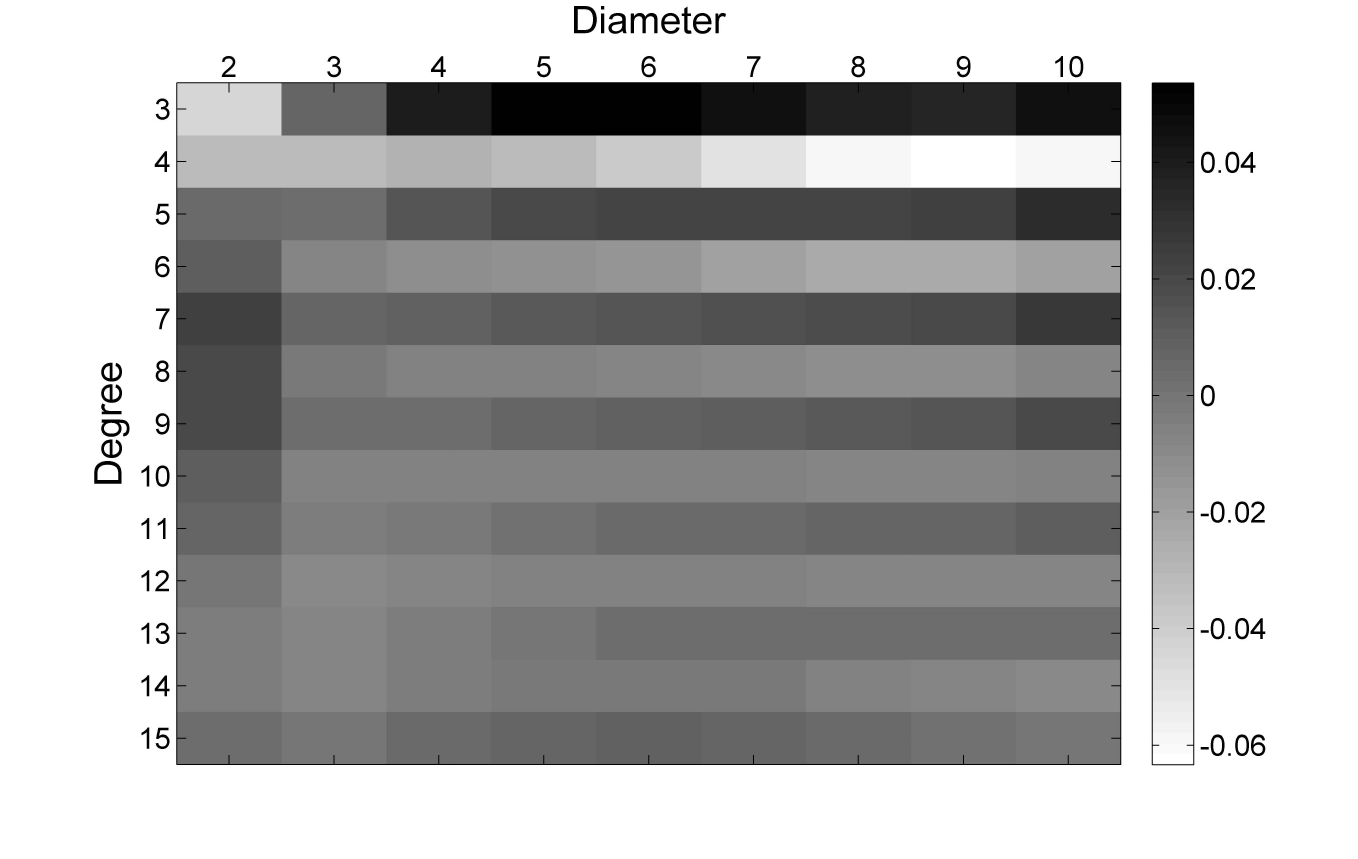}
\end{center}
\caption{Normalized difference between the predicted (logarithmic) upper bound minus actual (logarithmic) upper bound. Light entries correspond to negative differences. Darker colors correspond to positive differences.}
\label{fig:Log-upper-BW-1B-norma}
\end{figure}

If we normalize the difference at each $(\Delta, D)$ by the logarithm of the actual upper bound at $(\Delta, D)$, the disparity between even and odd rows is still visible, but it fades out as $\Delta$ increases, as shown in Figure \ref{fig:Log-upper-BW-1B-norma}. 

Circulant graphs of the form $C(n; \pm 1, \pm s_2, \ldots, \pm s_t)$ are called \emph{multi-loop graphs}. In particular, for $t = 2$ and $t = 3$ they are called \emph{double-loop graphs} and \emph{triple-loop graphs}, respectively. According to \cite{Mona12}, the maximum order of a triple-loop network $C(n; \pm 1, \pm s_2, \pm s_3)$ is:

\begin{equation}
\label{eq:triple}
T_{6,D} = \left\{ \begin{array}{ll}
\frac{32}{27}D^3 + \frac{16}{9}D^2 + 2D + 1 & \textrm{if $D \equiv 0 \bmod 3$}\\
32\lfloor \frac{D}{3} \rfloor^3 + 48\lfloor \frac{D}{3} \rfloor^2 + 30\lfloor \frac{D}{3} \rfloor + 7 & \textrm{if $D \equiv 1 \bmod 3$}\\
32\lfloor \frac{D}{3} \rfloor^3 + 80\lfloor \frac{D}{3} \rfloor^2 + 70\lfloor \frac{D}{3} \rfloor + 21 & \textrm{if $D \equiv 2 \bmod 3$}
\end{array} \right.
\end{equation}

Many of the graphs found by us here are multi-loop.

\section{A basic search procedure}
\label{sec:search}

Let $n$, $\Delta$ and $D$ be given. We now describe a systematic search procedure that attempts to find a circulant graph with the given parameters. Our search space consists of all connection sets $S$. The basic idea is to consider a tree-like organization of the search space, and then perform a Depth-First Search with backtracking on the tree, while pruning significant portions of it. To begin with, our search space will be restricted to connection sets with $s_1 = 1$. This will reduce the search space without losing too many potential solutions, and as a bonus, it will relieve us from the burden of testing connectedness. 

For even (resp. odd) $\Delta$, a node at depth $t$ in the search tree corresponds to a circulant graph of order $n$ and degree $2t$ (resp. $2t+1$). For the sake of brevity, we shall make no distinction between a node and its corresponding graph. The root of the search tree is the circulant graph $C(n; \pm1 \})$ (resp. $C(n;\{ \pm 1, n/2)$). Let $C(n;S)$ be a node of degree smaller than $\Delta$, and let $m$ be the largest element in $S$ such that $1 \le m < n/2$. Then $C(n; S)$ has offspring $C(n; S \cup \{\pm g\})$ for every $m+k \le g \le (n-k)/2$ (resp. $m+k \le g \le n/2-k$), where $k$ is a constant standing for the maximum difference allowed between two elements in $S$. Nodes of degree $\Delta$, appearing at depth $(\Delta-1)/2$ (resp. $(\Delta-3)/2$), are the leaves of the tree.

In the general Degree$/$Diameter Problem, it is well-known that large graphs tend to have large girth with respect to their diameter. This observation does not apply to circulant graphs, as their girth is always at most $4$ when $\Delta \ge 3$. Thus, we have generalized the girth criterion as follows:

For a circulant graph $C(n; S)$ and a length $d$, define $\PSet_w^d$ $(0 \le w < N)$ as the set of paths of length $d$ from a fixed vertex $v$ (say $0$) to the vertex $w$. Paths in $\PSet_w^d$ are meant to be disjoint, with no repeated vertices (except when $w = v$), and they are different up to commutation (that is, $P = Qx(xa)(xab)R$ and $P'=Qx(xb)(xba)R = Qx(xb)(xab)R$ are considered to be the same path).

Whether the offspring of a given node $C(n; S)$ are explored or not depends on a function on the cardinals of the sets $\PSet_w^d$. Basically, for each $i$ and $d$ we fix constants $c_i^d$ (for $0 \le i \le d < D$) which stand for the maximum number of paths allowed in $\PSet_w^d$ whenever $w$ is at distance $i \le d$ from the vertex $0$. That is, for every vertex $w$ in the graph, if $w$ is at distance $i < D$ from the vertex $0$, we require that $|\PSet_w^d| \le c_i^d$ for every $d$ with $i \le d < D$.

Tuning the constants $c_i^d$ provides great flexibility on the number of nodes explored, and proved to be very effective in practice. By defining small enough values for the constants $c_i^d$ it is possible to prune a significant portion of the search tree, especially as $\Delta$ increases. On the other hand, the constants $c_i^d$ should not be too small, so as not to miss too many potential solutions. 

Therefore, a careful choice of the constants $c_i^d$ is of paramount importance. In order to tune them properly, some degree of experimentation and intuition was involved. Initially, we tested many of the largest known circulant graphs for the values of $|\PSet_w^d|$. Most of these graphs were already known to be optimal.  

The parameter $k$, introduced to control the size of the search space, and thus the number of nodes visited by the algorithm, presented similar tradeoffs. Although certainly useful, experimenting with different values for the parameter $k$ provided rather poor improvements, compared to the cut-offs achievable via pruning (unless the ability to find existing solutions was seriously compromised by choosing $k$ to be too large).


Algorithm \ref{alg:search} formalizes the above ideas. Parameters $n$, $\Delta$ and $D$ are global. The input to the algorithm is the connection set $S$. Initially, the algorithm is called with $S = \{ \pm 1 \}$ (resp. $S = \{ \pm 1, n/2 \}$). Every time that a new generator is added, the algorithm calls itself recursively with the new connection set $S'$. If a circulant graph $C(n; S)$ with de desired parameters is found, the algorithm will print it out. Otherwise, the algorithm terminates without producing any output. 

Although the algorithm is presented in a recursive fashion, it was actually implemented using a stack for optimisation reasons. We omit actual implementation details here. 


\begin{algorithm}[H]
\SetKwInOut{Input}{Input}
\SetKwInOut{Output}{Output}
\vspace{.2cm}
\Input{A set of generators $S$.} 
\Output{A circulant graph $C(n; S)$ with degree $\Delta$ and diameter $D$.} 
\vspace{.2cm}
$m$:= Largest element in $S$\; 
\For{ $g$:= $m+k$ \KwTo $(n-k)/2$ } 
{
	 $S'$:= $S \cup \{\pm g\}$\;
	success := true\; 
	\For{each $g$-path $P$ in $C(n; S')$, with $P \in \PSet_w^d$}
	{
		$i$:= Distance from vertex $0$ to $w$ in $C(n; S')$\;
		\If{$d > i$ {\bf and} $|\PSet_w^d| > c_{i}^d$}
		{
			success := false\; 
			{\bf break}\;
		}
		\If{$d = i$ {\bf and} $|\PSet_w^j| > c_i^j$ for some $j$ ($i \le j < D$)}
		{
			success := false\; 
			{\bf break}\;
		}
	}
	\If{success}
	{
		\If{$|S'| < \Delta$}
		{
			Call SEARCH($S'$);
		}
		\Else
		{
			\If{diameter of $C(n; S')$ is $D$}{Print $C(n; S')$\;}
		}
	}
}

\caption{SEARCH}
\label{alg:search}  
\end{algorithm}

With the aid of Algorithm \ref{alg:search} we have been able to find $15$ new large circulant graphs of small diameter (ranging from $3$ to $5$), and degrees between $8$ and $15$. Additionally, we found two circulant graphs, which although not being record graphs, will be of interest later, since they will be used in the construction of larger graphs. The connection sets for all these graphs are listed in Table \ref{tab:ConnectionSets}. Four of them were found independently by Lewis in \cite{lewis14}, and proved to be optimal. These four graphs, as well as other graphs in Table \ref{tab:ConnectionSets}, up to degree 13, were also reported in \cite{enotes14}. By exhaustive search we also prove the optimality of the graph with degree $11$, diameter $3$, and order $210$, as well as the optimality of the previously known graphs of diameter $3$ and degrees $9$ and $10$.

\begin{table}[htp]
\begin{center}
\noindent\begin{tabular}{|r|r|lr|}
\hline

\rowcolor[gray]{.9}
$\Delta$ & $D$ & Order / Connection Set & \\
\hline

\rowcolor[gray]{1}
$8$ & $3$ & $104$ $\{\pm 1, \pm 16, \pm 20, \pm 27\}$ & $L,\bigstar$\\
\rowcolor[gray]{1}
       & $4$ & $248$ $\{\pm 1, \pm 61, \pm 72, \pm 76\}$ & $L,\bigstar$\\
\rowcolor[gray]{1}
       & $5$ & $528$ $\{\pm 1, \pm 89, \pm 156, \pm 162\}$ & $L,\bigstar$\\
\rowcolor[gray]{1}
       &        & $511$ $\{\pm 1, \pm 5, \pm 70, \pm 96\}$ & $\_,\bigstar$\\
\rowcolor[gray]{1}
       & $6$ & $967$ $\{\pm 1, \pm 7, \pm 132, \pm 182\}$ & $\_,\bigstar$\\

\rowcolor[gray]{.9}
$9$ & $4$ & $320$ $\{\pm 1, \pm 15, \pm 25, \pm 83, 160\}$ & $L,\bigstar$\\

\rowcolor[gray]{1}
$10$ & $4$ & $457$ $\{\pm 1, \pm 20, \pm 130, \pm 147, \pm 191\}$ & $\bigstar$\\
\rowcolor[gray]{1}
         & $5$ & $1099$ $\{\pm 1, \pm 53, \pm 207, \pm 272, \pm 536\}$ & $L$\\

\rowcolor[gray]{.9}
$11$ & $3$ & $210$ $\{\pm 1, \pm 49, \pm 59, \pm 84, \pm 89, 105\}$ & $\bigstar$\\
\rowcolor[gray]{.9}
         & $4$ & $576$ $\{\pm 1, \pm 9, \pm 75, \pm 155, \pm 179, 288\}$ & $\bigstar$\\
\rowcolor[gray]{.9}
         & $5$ & $1\,380$ $\{\pm 1, \pm 33, \pm 173, \pm 387, \pm 663, 690\}$ & $\bigstar$\\

\rowcolor[gray]{1}
$12$ & $3$ & $275$ $\{\pm 1, \pm 16, \pm 19, \pm 29, \pm 86, \pm 110\}$ & $\bigstar$\\
\rowcolor[gray]{1}
         & $4$ & $761$ $\{\pm 1, \pm 12, \pm 184, \pm 235, \pm 334, \pm 362\}$ & $\bigstar$\\
\rowcolor[gray]{1}
         & $5$ & $1\,800$ $\{\pm 1, \pm 30, \pm 64, \pm 384, \pm761, \pm 841\}$ & $\bigstar$\\

\rowcolor[gray]{.9}
$13$ & $3$ & $312$ $\{\pm 1, \pm 14, \pm 74, \pm 77, \pm 130, \pm 138, 156\}$ & $\bigstar$\\
\rowcolor[gray]{.9}
         & $4$ & $920$ $\{\pm 1, \pm 11, \pm 38, \pm 176, \pm 232, \pm 376, 460\}$ & $\bigstar$\\

\rowcolor[gray]{1}
$14$ & $3$ & $381$ $\{\pm 1, \pm 11, \pm 103, \pm 120, \pm 155, \pm 161, \pm 187\}$ & $\bigstar$\\

\rowcolor[gray]{.9}
$15$ & $3$ & $448$ $\{\pm 1, \pm 10, \pm 127, \pm 150, \pm 176, \pm 189, \pm 217, 224\}$ & $\bigstar$\\

\hline
\end{tabular}
\caption{Connection sets for graphs found by Algorithm \ref{alg:search}, labelled with a $\bigstar$. Graphs labelled with an $L$ were found independently by Lewis \cite{lewis14}, and proved to be optimal. Graphs labelled with a $\_$ symbol are sub-optimal. The graph with order $1099$ was also found by Lewis \cite{lewis14b}.}
\label{tab:ConnectionSets}
\end{center}
\end{table}

Unfortunately, for larger values of $\Delta$ and $D$, the execution time of Algorithm \ref{alg:search} becomes prohibitive with the computing power at our disposal. Thus, in order to complete our table of record circulant graphs up to $\Delta = 15$ and $D = 10$, we must use a combination of Algorithm \ref{alg:search} with Cartesian product, as described in the next section.

\section{Combining basic search with Cartesian product}
\label{sec:Cartesian}

Prior to our computer-based search, the main sources of large circulant graphs were:

\begin{enumerate}
\item For $t = 2$ an optimal circulant graph $C(n; \pm s_1, \pm s_2)$, is achieved for $s_1 = \lfloor \frac{1}{2}(\sqrt{2n-1}-1) \rfloor$ and $s_2 = s_1 + 1$ \cite{Ber89,BW85,Mona81}.
\item Monakhov and Monakhova used an evolutionary algorithm to find dense families of undirected circulant graphs. In particular, with the aid of this algorithm they found some families of large triple-loop graphs \cite{Mona02}. 
\item For larger degrees we have the construction $C(n; \pm 1, \pm s, \ldots, \pm s^{t-1})$, where $n = s^t$, and $s$ is an odd integer,  which yields good circulant graphs of diameter $\frac{t}{2}(s-1) = \frac{t}{2} n^{1/t} - \frac{t}{2}$ \cite{Wong74}.
\item Applying the methods described in \cite{Dou04} for abelian Cayley graphs, Charles Delorme has recently obtained several circulant graphs of large order \cite{Del13}.
\item With the aid of computer search, Lewis found several optimal graphs with degrees 8 and 9 \cite{lewis14}.
\item With the aid of Algorithm \ref{alg:search}, we found some new record graphs, up to degree 13 \cite{enotes14}. 
\end{enumerate}

\begin{figure}[htbp]
\begin{center}
	 	\includegraphics[width=0.6\textwidth]{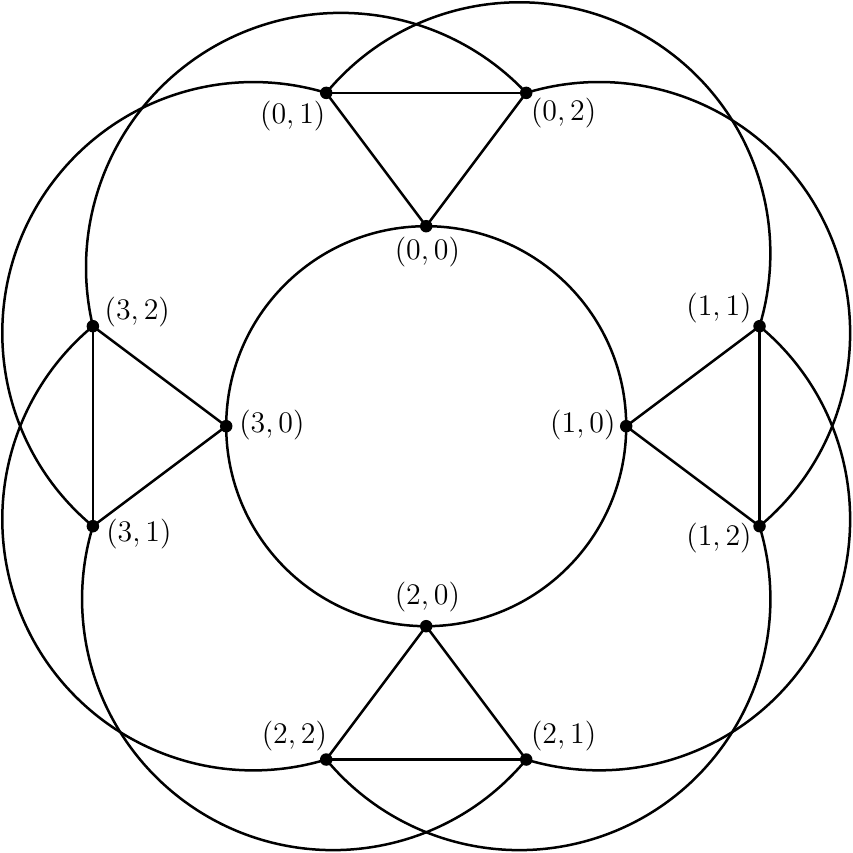}
\end{center}
\caption{Cartesian product $C(4; \pm 1) \ \Box \ C(3; \pm 1)$}
\label{fig:cartesian1}
\end{figure}

\begin{figure}[htbp]
\begin{center}
	 	\includegraphics[width=0.6\textwidth]{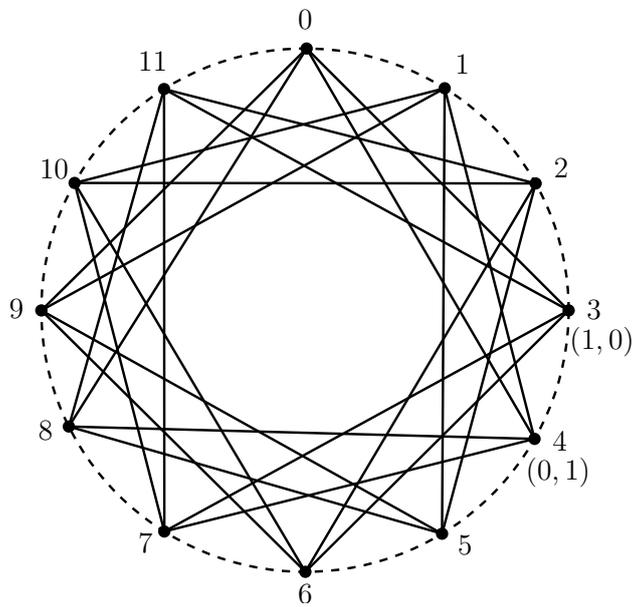}
\end{center}
\caption{The graph $C(12; \pm 3, \pm 4)$}
\label{fig:cartesian2}
\end{figure}

Our new computer search method is based on the combination of Algorithm \ref{alg:search} above with Cartesian product of graphs. If $n$ and $m$ are relatively prime, then $\mathbb{Z}_n \times \mathbb{Z}_m \cong \mathbb{Z}_{nm}$. Given the Cayley graphs $C(n; S_1)$ and $C(m; S_2)$, the Cayley graph $C(nm; S)$ can be obtained as $C(n; S_1) \ \Box \ C(m; S_2)$, where $S = m S_1 \cup n S_2$, and $\Box$ represents the Cartesian product of graphs. In this construction, the degree of $C(nm; S)$ is the sum of the degrees of the factor graphs $C(n; S_1)$ and $C(m; S_2)$, and the same holds for the diameter of $C(nm; S)$. Figure \ref{fig:cartesian1} depicts the Cartesian product $C(4; 1) \ \Box \ C(3; 1)$, and Figure \ref{fig:cartesian2} shows this same graph with a circulant layout. Note that the generator $3$ of  $\mathbb{Z}_{12}$ corresponds to $(1,0)$ in the --external-- direct product of $\mathbb{Z}_4$ and $\mathbb{Z}_3$, while the generator $4$ corresponds to $(0,1)$.


Now we just have to look for appropriate combinations of the largest known graphs in order to fill higher entries of the table. Normally we would be interested in using the largest possible factor graphs, but in order to apply the aforementioned construction we may need to resort to sub-optimal factor graphs, since the optimal factors may have orders which are not relatively prime, hence yielding a non-circulant graph. 

Algorithm \ref{alg:SearchAlgorithmCartesian} formalizes the method. The auxiliary function \emph{prime}($G$) determines whether $G$ is the Cartesian product of two circulant graphs whith relative prime orders. With the aid of Algorithm \ref{alg:SearchAlgorithmCartesian} we found $41$ new record circulant graphs, listed in Table \ref{tab:CartesianProduct}.

\clearpage

\begin{adjustbox}{max width=.90\textwidth}
\begin{algorithm}[H]
\SetKwInOut{Input}{Input}
\SetKwInOut{Output}{Output}
\vspace{.2cm}
\Input{A degree $\Delta \ge 4$ and a diameter $D \ge 10$.} 
\Output{Largest graph $C(nm; S) = C(n; S_1) \Box C(m; S_2)$ with degree $\Delta$ and diameter $D$; where $C(n; S_1)$ and $C(m; S_2)$ are  known circulant graphs, or are found by Algorithm \ref{alg:search}.} 
\vspace{.2cm}

$max$ := $0$\; 
$A$:= Empty array of graphs\;

\For{ $i$:= $2$ \KwTo $\lfloor \Delta/2 \rfloor$} 
{
	\For{ $j$:= $1$ \KwTo $\lfloor D/2 \rfloor$}
	{
		$G_1$ := Largest known circulant of deg. $i$ and diam. $j$\;
		$G_2$ := Largest known graph of deg. $\Delta - i$ and diam. $D - j$\;
		
		\If{$|G_1||G_2| > max$}
		{
			Insert $|G_1| \Box |G_2|$ in $A$ so as to keep $A$ sorted\;

			\If{$gcd(|G_1|,|G_2|) = 1$}
			{
				$max$ := $|G_1||G_2|$\; 
				Remove from $A$ all the elements after $|G_1| \Box |G_2|$\;
			}
		}
	} 
}

\While{$length(A) > 0$ {\bf and not} $prime(A[0])$}
{
	$G$ := $|G_1| \Box |G_2|$ := $A[0]$\;
	Remove $G$ from $A$\;

	$n'$ := $|G_1|-1$\;
	\While{$n'|G_2| > max$}
	{
		\If{$gcd(n',m) = 1$}
		{
			$G'_1$ := Circulant of order $n'$, same deg. and diam. as $G_1$\; 
			{\bf if } $G'_1$ was found {\bf then break}\;
			$n'$ := $n'-1$\;
		}
	}

	$m'$ := $|G_2|-1$\;
	\While{$|G_1|m' > max$}
	{
		\If{$gcd(n,m') = 1$}
		{
			$G'_2$ := Circulant of order $m'$, same deg. and diam. as $G_2$\; 
			{\bf if } $G'_2$ was found {\bf then break}\;
			$m'$ := $m'-1$\;
		}
	}

	$G'$ := $n'm > nm'$ {\bf ?} $G'_1 \Box G_2$ {\bf :} $G_1 \Box G'_2$\;
	Insert $G'$ in $A$ so as to keep $A$ sorted\;
	$max$ := $|G'|$\;
	Remove from $A$ all the elements after $G'$\;
}

{\bf if } $length(A) > 0$ {\bf then} Print $A[0]$\;

\caption{COMBINED SEARCH AND CARTESIAN PRODUCT}
\label{alg:SearchAlgorithmCartesian}  
\end{algorithm}
\end{adjustbox}

\begin{table}[htp]
\begin{center}
\noindent\begin{tabular}{|rrr|rrr|rrr|}
\hline

\rowcolor[gray]{.9}
\multicolumn{3}{|c|}{$G_1 \Box G_2$} & \multicolumn{3}{c|}{$G_1$} & \multicolumn{3}{c|}{$G_2$} \\
\rowcolor[gray]{.9}
$\Delta$ & $D$ & Order & $\Delta$ & $D$ & Order & $\Delta$ & $D$ & Order\\
\hline

\rowcolor[gray]{1}
$10$ & $6$ & ${\bf 1\,533}$    &    $2$ & $1$ & ${\bf 3}$    &    $8$ & $5$ & ${\bf 511}$\\
\rowcolor[gray]{1}
 	& $7$ & ${\bf 2\,925}$    &    $4$ & $3$ & ${\bf 25}$    &    $6$ & $4$ & ${\bf 117}$\\
\rowcolor[gray]{1}
 	& $8$ & ${\bf 5\,136}$    &    $2$ & $1$ & ${\bf 3}$    &    $8$ & $7$ & ${\bf 1\,712}$\\
\rowcolor[gray]{1}
 	& $9$ & ${\bf 8\,560}$    &    $2$ & $2$ & ${\bf 5}$    &    $8$ & $7$ & ${\bf 1\,712}$\\
\rowcolor[gray]{1}
	& $10$ & ${\bf 13\,840}$    &    $2$ & $2$ & ${\bf 5}$    &    $8$ & $8$ & ${\bf 2\,768}$\\

\rowcolor[gray]{.9}
$11$ & $6$ & ${\bf 2\,100}$    &    $2$ & $1$ & ${\bf 3}$    &    $9$ & $5$ & ${\bf 700}$\\
\rowcolor[gray]{.9}
 	& $7$ & ${\bf 4\,088}$    &    $3$ & $2$ & ${\bf 8}$    &    $8$ & $5$ & ${\bf 511}$\\
\rowcolor[gray]{.9}
 	& $8$ & ${\bf 7\,736}$    &    $3$ & $2$ & ${\bf 8}$    &    $8$ & $6$ & ${\bf 984}$\\
\rowcolor[gray]{.9}
 	& $9$ & ${\bf 13\,400}$    &    $4$ & $3$ & ${\bf 25}$    &    $7$ & $6$ & ${\bf 536}$\\
\rowcolor[gray]{.9}
	& $10$ & ${\bf 21\,976}$    &    $4$ & $4$ & ${\bf 41}$    &    $7$ & $6$ & ${\bf 536}$\\

\rowcolor[gray]{1}
$12$ & $6$ & ${\bf 3\,297}$    &    $2$ & $1$ & ${\bf 3}$    &    $10$ & $5$ & ${\bf 1\,099}$\\
\rowcolor[gray]{1}
 	& $7$ & ${\bf 6\,864}$    &    $4$ & $2$ & ${\bf 13}$    &    $8$ & $5$ & ${\bf 528}$\\
\rowcolor[gray]{1}
 	& $8$ & ${\bf 13\,200}$    &    $4$ & $3$ & ${\bf 25}$    &    $8$ & $5$ & ${\bf 528}$\\
\rowcolor[gray]{1}
 	& $9$ & ${\bf 24\,600}$    &    $4$ & $3$ & ${\bf 25}$    &    $8$ & $6$ & ${\bf 984}$\\
\rowcolor[gray]{1}
	& $10$ & ${\bf 42\,800}$    &    $4$ & $3$ & ${\bf 25}$    &    $8$ & $7$ & ${\bf 1\,712}$\\

\rowcolor[gray]{.9}
$13$ & $5$ & ${\bf 1\,828}$    &    $3$ & $1$ & ${\bf 4}$    &    $10$ & $4$ & ${\bf 457}$\\
\rowcolor[gray]{.9}
	& $6$ & ${\bf 4\,396}$    &    $3$ & $1$ & ${\bf 4}$    &    $10$ & $5$ & ${\bf 1\,099}$\\
\rowcolor[gray]{.9}
 	& $7$ & ${\bf 9\,100}$    &    $4$ & $2$ & ${\bf 13}$    &    $9$ & $5$ & ${\bf 700}$\\
\rowcolor[gray]{.9}
 	& $8$ & ${\bf 18\,720}$    &    $6$ & $4$ & ${\bf 117}$    &    $7$ & $4$ & ${\bf 160}$\\
\rowcolor[gray]{.9}
 	& $9$ & ${\bf 36\,036}$    &    $6$ & $4$ & ${\bf 117}$    &    $7$ & $5$ & ${\bf 308}$\\
\rowcolor[gray]{.9}
	& $10$ & ${\bf 63\,700}$    &    $4$ & $3$ & ${\bf 25}$    &    $9$ & $7$ & ${\bf 2\,548}$\\

\rowcolor[gray]{1}
$14$ & $4$ & ${\bf 825}$    &    $2$ & $1$ & ${\bf 3}$    &    $12$ & $3$ & ${\bf 275}$\\
\rowcolor[gray]{1}
 	& $5$ & ${\bf 2\,285}$    &    $4$ & $1$ & ${\bf 5}$    &    $10$ & $4$ & ${\bf 457}$\\
\rowcolor[gray]{1}
	& $6$ & ${\bf 5\,941}$    &    $4$ & $2$ & ${\bf 13}$    &    $10$ & $4$ & ${\bf 457}$\\
\rowcolor[gray]{1}
 	& $7$ & ${\bf 14\,287}$    &    $4$ & $2$ & ${\bf 13}$    &    $10$ & $5$ & ${\bf 1\,099}$\\
\rowcolor[gray]{1}
 	& $8$ & ${\bf 29\,016}$    &    $6$ & $4$ & ${\bf 117}$    &    $8$ & $4$ & ${\bf 248}$\\
\rowcolor[gray]{1}
 	& $9$ & ${\bf 59\,787}$    &    $6$ & $4$ & ${\bf 117}$    &    $8$ & $5$ & ${\bf 511}$\\
\rowcolor[gray]{1}
	& $10$ & ${\bf 113\,139}$    &    $6$ & $4$ & ${\bf 117}$    &    $8$ & $6$ & ${\bf 967}$\\

\rowcolor[gray]{.9}
$15$ & $4$ & ${\bf 1\,100}$    &    $3$ & $1$ & ${\bf 4}$    &    $12$ & $3$ & ${\bf 275}$\\
\rowcolor[gray]{.9}
	& $5$ & ${\bf 3\,044}$    &    $3$ & $1$ & ${\bf 4}$    &    $12$ & $4$ & ${\bf 761}$\\
\rowcolor[gray]{.9}
	& $6$ & ${\bf 7\,524}$    &    $7$ & $3$ & ${\bf 76}$    &    $8$ & $3$ & ${\bf 99}$\\
\rowcolor[gray]{.9}
 	& $7$ & ${\bf 17\,940}$    &    $4$ & $2$ & ${\bf 13}$    &    $11$ & $5$ & ${\bf 1\,380}$\\
\rowcolor[gray]{.9}
 	& $8$ & ${\bf 39\,564}$    &    $5$ & $3$ & ${\bf 36}$    &    $10$ & $5$ & ${\bf 1\,099}$\\
\rowcolor[gray]{.9}
 	& $9$ & ${\bf 81\,900}$    &    $6$ & $4$ & ${\bf 117}$    &    $9$ & $5$ & ${\bf 700}$\\
\rowcolor[gray]{.9}
	& $10$ & ${\bf 154\,720}$    &    $7$ & $4$ & ${\bf 160}$    &    $8$ & $6$ & ${\bf 967}$\\

\rowcolor[gray]{1}
$16$ & $5$ & ${\bf 3\,805}$    &    $4$ & $1$ & ${\bf 5}$    &    $12$ & $4$ & ${\bf 761}$\\
\rowcolor[gray]{1}
	& $6$ & ${\bf 10\,296}$    &    $8$ & $3$ & ${\bf 99}$    &    $8$ & $3$ & ${\bf 104}$\\
\rowcolor[gray]{1}
 	& $7$ & ${\bf 25\,135}$    &    $6$ & $3$ & ${\bf 55}$    &    $10$ & $4$ & ${\bf 457}$\\
\rowcolor[gray]{1}
 	& $8$ & ${\bf 60\,445}$    &    $6$ & $3$ & ${\bf 55}$    &    $10$ & $5$ & ${\bf 1\,099}$\\
\rowcolor[gray]{1}
 	& $9$ & ${\bf 128\,583}$    &    $6$ & $4$ & ${\bf 117}$    &    $10$ & $5$ & ${\bf 1\,099}$\\
\rowcolor[gray]{1}
	& $10$ & ${\bf 269\,808}$    &    $8$ & $5$ & ${\bf 511}$    &    $8$ & $5$ & ${\bf 528}$\\

\hline
\end{tabular}
\caption{Large circulant graphs obtained by the Cartesian product of two circulant graphs. 
}
\label{tab:CartesianProduct}
\end{center}
\end{table}

\section{The largest known circulant graphs}
\label{sec:table}

Table \ref{tab:circulant1} summarizes our current knowledge about the largest known circulant graphs, including the ones obtained with the aid of our methods. Each entry of the table contains four elements: the order of the largest circulant graph known to-date (upper left), a reference to its source (upper right), the best-known upper bound (lower left), and the percentage of the upper bound achieved (lower right). This table can also be found in \cite{circulant_tables}.

\begin{sidewaystable}
\scalebox{0.85}{ 
\begin{tabular}{|cc|*{9}{ll|}} \hline 
\multicolumn{2}{|c|}{} & \multicolumn{18}{c|}{\large{Diameter $D$}} \\
\rowcolor[gray]{.9} 
\multicolumn{2}{|c|}{} & \multicolumn{2}{c}{{\bf 2}} & \multicolumn{2}{c}{{\bf 3}} & \multicolumn{2}{c}{{\bf 4}} & \multicolumn{2}{c}{{\bf 5}} & \multicolumn{2}{c}{{\bf 6}} & \multicolumn{2}{c}{{\bf 7}} & \multicolumn{2}{c}{{\bf 8}} & \multicolumn{2}{c}{{\bf 9}} & \multicolumn{2}{c|}{{\bf 10}} \\ \hline
 & {\bf 3} & 8 & & 12 & & 16 & & 20 & & 24 & & 28 & & 32 & & 36 & & 40 & \\
 & & 8 & 100\% & 12 & 100\% & 16 & 100\% & 20 & 100\% & 24 & 100\% & 28 & 100\% & 32 & 100\% & 36 & 100\% & 40 & 100\% \\
\rowcolor[gray]{.9} 
 & {\bf 4} & 13 & & 25 & & 41 & & 61 & & 85 & & 113 & & 145 & & 181 & & 221 & \\
\rowcolor[gray]{.9}
 & & 13 & 100\% & 25 & 100\% & 41 & 100\% & 61 & 100\% & 85 & 100\% & 113 & 100\% & 145 & 100\% & 181 & 100\% & 221 & 100\% \\ 
 & {\bf 5} & 16 & & 36 & \cite{Mac10} & 64 & \cite{Mac10} & 100 & \cite{Mac10} & 144 & \cite{Mac10} & 196 & \cite{Mac10} & 256 & \cite{Mac10} & 324 & \cite{Mac10} & 400 & \cite{Mac10} \\
 & & 18 & 89\% & 38 & 95\% & 66 & 97\% & 102 & 98\% & 146 & 99\% & 198 & 99\% & 258 & 99\% & 326 & 99\% & 402 & 99\% \\  
\rowcolor[gray]{.9}
 & {\bf 6} & 21 & \cite{Del13} & 55 & \cite{Del13} & 117 & \cite{Del13} & 203 & \cite{Del13} & 333 & \cite{Mac10} & 515 & \cite{Mac10} & 737 & \cite{Mac10} & 1027 & \cite{Mac10} & 1393 & \cite{Mona12} \\
\rowcolor[gray]{.9}
 & & 25 & 84\% & 63 & 87\% & 129 & 91\% & 231 & 88\% & 377 & 88\% & 575 & 90\% & 833 & 88\% & 1159 & 89\% & 1561 & 89\% \\  
 & {\bf 7} & 26 & \cite{Del13} & 76 & \cite{Mac10} & 160 & \cite{Mac10} & 308 & \cite{Mac10} & 536 & T.\ref{tab:circulant2} & 828 & T.\ref{tab:circulant2} & 1232 & T.\ref{tab:circulant2} & 1764 & T.\ref{tab:circulant2} & 2392 & T.\ref{tab:circulant2} \\
 & & 32 & 81\% & 88 & 86\% & 192 & 83\% & 360 & 86\% & 608 & 88\% & 952 & 87\% & 1408 & 87\% & 1992 & 88\% & 2720 & 88\% \\  
\rowcolor[gray]{.9}
 & {\bf 8} & 35 & \cite{Del13} & 104 & T.\ref{tab:circulant2} & 248 & T.\ref{tab:circulant2} & 528 & T.\ref{tab:circulant2} & 984 & \cite{lewis14} & 1712 & \cite{lewis14} & 2768 & \cite{lewis14} & 4280 & \cite{lewis14} & 6320 & \cite{lewis14} \\
\rowcolor[gray]{.9}
 & & 41 & 85\% & 129 & 77\% & 321 & 75\% & 681 & 77\% & 1289 & 76\% & 2241 & 76\% & 3649 & 76\% & 5641 & 76\% & 8361 & 76\% \\  
\Large{$\Delta$} & {\bf 9} & 42 & \cite{Del13} & 130 & \cite{Mac10} & 320 & T.\ref{tab:circulant2} & 700 & \cite{lewis14} & 1416 & \cite{lewis14} & 2548 & \cite{lewis14} & 4304 & \cite{lewis14} & 6804 & \cite{lewis14} & 10320 & \cite{lewis14} \\
 & & 50 & 84\% & 170 & 76\% & 450 & 71\% & 1002 & 70\% & 1970 & 72\% & 3530 & 72\% & 5890 & 73\% & 9290 & 73\% & 14002 & 74\% \\
\rowcolor[gray]{.9}
 & {\bf 10} & 51 & \cite{Del13} & 177 & \cite{Mac10} & 457 & T.\ref{tab:circulant2} & 1099 & \cite{lewis14b} & 1533 & T.\ref{tab:circulant2} & 2925 & T.\ref{tab:circulant2} & 5136 & T.\ref{tab:circulant2} & 8560 & T.\ref{tab:circulant2} & 13840 & T.\ref{tab:circulant2} \\
\rowcolor[gray]{.9}
 & & 61 & 84\% & 231 & 77\% & 681 & 67\% & 1683 & 65\% & 3653 & 42\% & 7183 & 40\% & 13073 & 39\% & 22363 & 38\% & 36365 & 38\% \\ 
 & {\bf 11} & 56 & \cite{Del13} & 210 & T.\ref{tab:circulant2} & 576 & T.\ref{tab:circulant2} & 1380 & T.\ref{tab:circulant2} & 2100 & T.\ref{tab:circulant2} & 4088 & T.\ref{tab:circulant2} & 7736 & T.\ref{tab:circulant2} & 13400 & T.\ref{tab:circulant2} & 21976 & T.\ref{tab:circulant2} \\
 & & 72 & 78\% & 292 & 72\% & 912 & 63\% & 2364 & 58\% & 5336 & 39\% & 10836 & 37\% & 20256 & 38\% & 35436 & 38\% & 58728 & 37\% \\  
\rowcolor[gray]{.9}
 & {\bf 12} & 67 & \cite{Del13} & 275 & T.\ref{tab:circulant2} & 761 & T.\ref{tab:circulant2} & 1800 & T.\ref{tab:circulant2} & 3297 & T.\ref{tab:circulant2} & 6864 & T.\ref{tab:circulant2} & 13200 & T.\ref{tab:circulant2} & 24600 & T.\ref{tab:circulant2} & 42800 & T.\ref{tab:circulant2} \\ 
\rowcolor[gray]{.9} 
 & & 85 & 79\% & 377 & 73\% & 1289 & 59\% & 3653 & 49\% & 8989 & 36\% & 19825 & 34\% & 40081 & 33\% & 75517 & 32\% & 134245 & 32\% \\ 
 & {\bf 13} & 80 & T.\ref{tab:circulant3} & 312 & T.\ref{tab:circulant3} & 920 & T.\ref{tab:circulant3} & 1828 & T.\ref{tab:circulant3} & 4396 & T.\ref{tab:circulant3} & 9100 & T.\ref{tab:circulant3} & 18720 & T.\ref{tab:circulant3} & 36036 & T.\ref{tab:circulant3} & 63700 & T.\ref{tab:circulant3} \\
 & & 98 & 82\% & 462 & 68\% & 1666 & 55\% & 4942 & 37\% & 12642 & 35\% & 28814 & 31\% & 59906 & 31\% & 115598 & 31\% & 209762 & 30\% \\ 
\rowcolor[gray]{.9}
 & {\bf 14} & 90 & T.\ref{tab:circulant3} & 381 & T.\ref{tab:circulant3} & 825 & T.\ref{tab:circulant3} & 2285 & T.\ref{tab:circulant3} & 5941 & T.\ref{tab:circulant3} & 14287 & T.\ref{tab:circulant3} & 29016 & T.\ref{tab:circulant3} & 54120 & T.\ref{tab:circulant3} & 113139 & T.\ref{tab:circulant3} \\ 
\rowcolor[gray]{.9}
 & & 113 & 80\% & 575 & 66\% & 2241 & 37\% & 7183 & 32\% & 19825 & 30\% & 48639 & 30\% & 108545 & 27\% & 224143 & 24\% & 433905 & 26\% \\ 
 & {\bf 15} & 96 & T.\ref{tab:circulant3} & 448 & T.\ref{tab:circulant3} & 1100 & T.\ref{tab:circulant3} & 2880 & T.\ref{tab:circulant3} & 7488 & T.\ref{tab:circulant3} & 17584 & T.\ref{tab:circulant3} & 39564 & T.\ref{tab:circulant3} & 81900 & T.\ref{tab:circulant3} & 154720 & T.\ref{tab:circulant3} \\
 & & 128 & 75\% & 688 & 65\% & 2816 & 39\% & 9424 & 30\% & 27008 & 28\% & 68464 & 25\% & 157184 & 25\% & 332688 & 24\% & 658048 & 23\% \\ 
 \rowcolor[gray]{.9}
 & {\bf 16} & 112 & T.\ref{tab:circulant2} & xxx & T.\ref{tab:circulant3} & 936 & T.\ref{tab:circulant3} & 3640 & T.\ref{tab:circulant3} & 9597 & T.\ref{tab:circulant3} & 25135 & T.\ref{tab:circulant3} & 60445 & T.\ref{tab:circulant3} & 128583 & T.\ref{tab:circulant3} & 239816 & T.\ref{tab:circulant3} \\ 
 \rowcolor[gray]{.9}
  & & 145 & 25\% & 833 & 11\% & 3649 & 25\% & 13073 & 28\% & 40081 & 24\% & 108545 & 23\% & 265729 & 23\% & 598417 & 21\% & 1256465 & 19\% \\ \hline
\end{tabular} }
\caption{Orders of the largest known circulant graphs.} 
\label{tab:circulant1}
\end{sidewaystable}

The actual graphs obtained by Algorithms \ref{alg:search} and \ref{alg:SearchAlgorithmCartesian} are given in tables \ref{tab:circulant2} and \ref{tab:circulant3}. The tables also indicate which of these graphs are optimal multi-loop graphs. Optimality has been verified for some small values of $\Delta$ and $D$, by performing a systematic  search through all multi-loop graphs with the given parameters $\Delta$ and $D$.   

\begin{table}[htp]
\begin{center}
\scalebox{0.9}{ 
\begin{tabular}{|lll|} \hline
$\Delta$ & $D$ & Circulant graph \\ \hline \hline 
\rowcolor[gray]{.9} 
7 & 6 & $536\{\pm 1, \pm 231, \pm 239, 268 \} \; \star$ \\ 
7 & 7 & $828\{\pm 1, \pm 9, \pm 91, 414 \} \; \star$ \\ 
\rowcolor[gray]{.9} 
7 & 8 & $1232\{\pm 1, \pm 11, \pm 111, 616 \} \; \star$ \\ 
7 & 9 & $1764\{\pm 1, \pm 803, \pm 815, 882 \} \; \star$ \\ 
\rowcolor[gray]{.9} 
7 & 10 & $2392\{\pm 1, \pm 13, \pm 183, 1196 \} \; \star$ \\ 
8 & 3 & $104\{\pm 1, \pm 16, \pm 20, \pm 27 \} \; \star$ \\ 
\rowcolor[gray]{.9} 
8 & 4 & $248\{\pm 1, \pm 61, \pm 72, \pm 76 \} \; \star$ \\ 
8 & 5 & $528\{\pm 1, \pm 89, \pm 156, \pm 162 \}$ \\ 
\rowcolor[gray]{.9} 
8 & 6 & $967\{\pm 1, \pm 7, \pm 132, \pm 182 \}$ \\ 
8 & 7 & $1545\{\pm 1, \pm 170, \pm 178, \pm 468 \}$ \\ 
\rowcolor[gray]{.9} 
9 & 4 & $320\{\pm 1, \pm 15, \pm 25, \pm 83, 160 \} \; \star$ \\ 
9 & 5 & $684\{\pm 1, \pm 111, \pm 145, \pm 279, 342 \}$ \\ 
\rowcolor[gray]{.9} 
9 & 6 & $1284\{\pm 1, \pm 36, \pm 163, \pm 342, 642 \}$ \\ 
9 & 7 & $2340\{\pm 1, \pm 149, \pm 157, \pm 645, 1170 \}$ \\ 
\rowcolor[gray]{.9} 
10 & 4 & $457\{\pm 1, \pm 20, \pm 130, \pm 147, \pm 191 \}$ \\ 
10 & 5 & $1099\{\pm 1, \pm 53, \pm 207, \pm 272, \pm 536 \}$ \\ 
\rowcolor[gray]{.9}
10 & 6 & $1533\{\pm 3, \pm 15, \pm 210, \pm 288, \pm 511 \}$ \\ 
10 & 7 & $2925\{\pm 25, \pm 351, \pm 400, \pm 468, \pm 550 \}$ \\ 
\rowcolor[gray]{.9}
10 & 8 & $5136\{\pm 3, \pm 645, \pm 1712, \pm 1824, \pm 1848 \}$ \\ 
10 & 9 & $8560\{\pm 5, \pm 1075, \pm 1712, \pm 3040, \pm 3080 \}$ \\ 
\rowcolor[gray]{.9}
10 & 10 & $13840\{\pm 5, \pm 1032, \pm 2768, \pm 5360, \pm 5400 \}$ \\ 
11 & 3 & $210\{\pm 1, \pm 49, \pm 59, \pm 84, \pm 89, 105 \} \; \star$ \\ 
\rowcolor[gray]{.9} 
11 & 4 & $576\{\pm 1, \pm 9, \pm 75, \pm 155, \pm 179, 288 \}$ \\ 
11 & 5 & $1380\{\pm 1, \pm 33, \pm 173, \pm 387, \pm 663, 690 \}$ \\ 
\rowcolor[gray]{.9} 
11 & 6 & $2100\{\pm 3, \pm 15, \pm 591, \pm 669, \pm 700, 1050 \}$ \\ 
11 & 7 & $4088\{\pm 8, \pm 40, \pm 511, \pm 560, \pm 768, 2044 \}$ \\ 
\rowcolor[gray]{.9}
11 & 8 & $7736\{\pm 8, \pm 56, \pm 967, \pm 1056, \pm 1456, 3868 \}$ \\ 
11 & 9 & $13400\{\pm 25, \pm 1608, \pm 2144, \pm 5775, \pm 5975, 6700 \}$ \\ 
\rowcolor[gray]{.9} 
11 & 10 & $21976\{\pm 41, \pm 2144, \pm 2680, \pm 9471, \pm 9799, 10988 \}$ \\ 
12 & 3 & $275\{\pm 1, \pm 16, \pm 19, \pm 29, \pm 86, \pm 110 \}$ \\ 
\rowcolor[gray]{.9}
12 & 4 & $761\{\pm 1, \pm 12, \pm 184, \pm 235, \pm 334, \pm 362 \}$ \\ 
12 & 5 & $1800\{\pm 1, \pm 30, \pm 64, \pm 384, \pm 761, \pm 841 \}$ \\ 
\rowcolor[gray]{.9}
12 & 6 & $3297\{\pm 3, \pm 159, \pm 621, \pm 816, \pm 1099, \pm 1608 \}$ \\ 
12 & 7 & $6864\{\pm 13, \pm 1056, \pm 1157, \pm 1584, \pm 2028, \pm 2106 \}$ \\ 
\rowcolor[gray]{.9}
12 & 8 & $13200\{\pm 25, \pm 1584, \pm 2112, \pm 2225, \pm 3900, \pm 4050 \}$ \\ 
12 & 9 & $24600\{\pm 25, \pm 2952, \pm 3936, \pm 4075, \pm 8700, \pm 8850 \}$ \\ 
\rowcolor[gray]{.9}
12 & 10 & $42800\{\pm 25, \pm 5136, \pm 5375, \pm 6848, \pm 15200, \pm 15400 \}$ \\ 
\hline 
\end{tabular} } 
\caption{New circulant networks obtained by Algorithms \ref{alg:search} and \ref{alg:SearchAlgorithmCartesian}. The graphs marked with $\star$ are optimal.}
\label{tab:circulant2}
\end{center}
\end{table}

\begin{table}[htp]
\begin{center}
\scalebox{0.9}{ 
\begin{tabular}{|lll|} \hline
$\Delta$ & $D$ & Circulant graph \\ \hline \hline 
\rowcolor[gray]{.9} 
13 & 2 & $80\{\pm 1, \pm 3, \pm 9, \pm 20, \pm 25, \pm 33, 40 \}$ \\ 
13 & 3 & $312\{\pm 1, \pm 14, \pm 74, \pm 77, \pm 130, \pm 138, 156 \}$ \\ 
\rowcolor[gray]{.9}
13 & 4 & $920\{\pm 1, \pm 11, \pm 38, \pm 176, \pm 232, \pm 376, 460 \}$ \\ 
13 & 5 & $1828\{\pm 4, \pm 80, \pm 457, \pm 520, \pm 588, \pm 764, 914 \}$ \\ 
\rowcolor[gray]{.9}
13 & 6 & $4396\{\pm 4, \pm 212, \pm 828, \pm 1088, \pm 1099, \pm 2144, 2198 \}$ \\ 
13 & 7 & $9100\{\pm 13, \pm 65, \pm 1400, \pm 2100, \pm 2561, \pm 2899, 4550 \}$ \\ 
\rowcolor[gray]{.9}
13 & 8 & $18720\{\pm 117, \pm 160, \pm 585, \pm 2560, \pm 3520, \pm 3627, 9360 \}$ \\ 
13 & 9 & $36036\{\pm 117, \pm 308, \pm 819, \pm 4928, \pm 5031, \pm 6776, 18018 \}$ \\ 
\rowcolor[gray]{.9} 
13 & 10 & $63700\{\pm 25, \pm 175, \pm 7644, \pm 10192, \pm 13025, \pm 14275, 31850 \}$ \\ 
14 & 2 & $90\{\pm 1, \pm 4, \pm 10, \pm 17, \pm 26, \pm 29, \pm 41 \} \; \star$ \\ 
\rowcolor[gray]{.9} 
14 & 3 & $381\{\pm 1, \pm 11, \pm 103, \pm 120, \pm 155, \pm 161, \pm 187 \} \; \star$ \\ 
14 & 4 & $825\{\pm 3, \pm 48, \pm 57, \pm 87, \pm 258, \pm 275, \pm 330 \} \; \star$ \\ 
\rowcolor[gray]{.9}
14 & 5 & $2285\{\pm 5, \pm 100, \pm 457, \pm 650, \pm 735, \pm 914, \pm 955 \} \; \star$ \\ 
14 & 6 & $5941\{\pm 13, \pm 260, \pm 914, \pm 1371, \pm 1690, \pm 1911, \pm 2483 \} \; \star$ \\ 
\rowcolor[gray]{.9} 
14 & 7 & $14287\{\pm 13, \pm 689, \pm 2198, \pm 2691, \pm 3297, \pm 3536, \pm 6968, \} \; \star$ \\ 
14 & 8 & $29016\{\pm 117, \pm 248, \pm 3968, \pm 5456, \pm 7137, \pm 8424, \pm 8892 \} \; \star$ \\ 
\rowcolor[gray]{.9} 
14 & 9 & $54120\{\pm 55, \pm 984, \pm 4920, \pm 8965, \pm 19140, \pm 19470, \pm 20664 \}$ \\ 
14 & 10 & $113139\{\pm 117, \pm 819, \pm 967, \pm 15444, \pm 15472, \pm 21274, \pm 21294 \}$ \\ 
\rowcolor[gray]{.9}
15 & 2 & $96\{\pm 1, \pm 3, \pm 5, \pm 11, \pm 24, \pm 31, \pm 39, 48 \}$ \\ 
15 & 3 & $448\{\pm 1, \pm 10, \pm 127, \pm 150, \pm 176, \pm 189, \pm 217, 224 \} \; \star$ \\ 
\rowcolor[gray]{.9} 
15 & 4 & $1100\{\pm 4, \pm 64, \pm 76, \pm 116, \pm 275, \pm 344, \pm 440, 550 \}$ \\ 
15 & 5 & $2880\{\pm 5, \pm 45, \pm 375, \pm 576, \pm 775, \pm 895, \pm 1152, 1440 \}$ \\ 
\rowcolor[gray]{.9} 
15 & 6 & $7488\{\pm 13, \pm 117, \pm 975, \pm 1152, \pm 1728, \pm 2015, \pm 2327, 3744 \}$ \\ 
15 & 7 & $17584\{\pm 16, \pm 848, \pm 1099, \pm 3297, \pm 3312, \pm 4352, \pm 8576, 8792 \}$ \\ 
\rowcolor[gray]{.9}
15 & 8 & $39564\{\pm 36, \pm 1099, \pm 1908, \pm 5495, \pm 7452, \pm 9792, \pm 19296, 19782 \}$ \\ 
15 & 9 & $81900\{\pm 117, \pm 585, \pm 700, \pm 11200, \pm 15400, \pm 23049, \pm 26091, 40950 \}$ \\ 
\rowcolor[gray]{.9}
15 & 10 & $154720\{\pm 160, \pm 967, \pm 1120, \pm 4835, \pm 21120, \pm 29120, \pm 29977, 77360 \}$ \\ 
\hline 
\end{tabular} } 
\caption{New circulant networks obtained by Algorithms \ref{alg:search} and \ref{alg:SearchAlgorithmCartesian}. The graphs marked with $\star$ are optimal.}
\label{tab:circulant3}
\end{center}
\end{table}

Table \ref{tab:circulant1} can be well approximated by a bi-exponential function, such as $e^{f(\Delta, D)}$, where 

\begin{equation}
\label{eq:fit2}
\begin{split}
f(\Delta, D) = &-1.912 + 0.8845 \Delta + 0.8014 D -0.08641 \Delta^2 \\
               &+ 0.1152 \Delta D -0.1214 D^2 + 0.002799 \Delta^3 \\
               &+ 0.002799 \Delta^2 D -0.001046 \Delta D^2 + 0.005036 D^3,
\end{split}
\end{equation}

with $R^2 = 0.997$. 

\begin{figure}[htbp]
\begin{center}
	 	\includegraphics[width=1.0\textwidth]{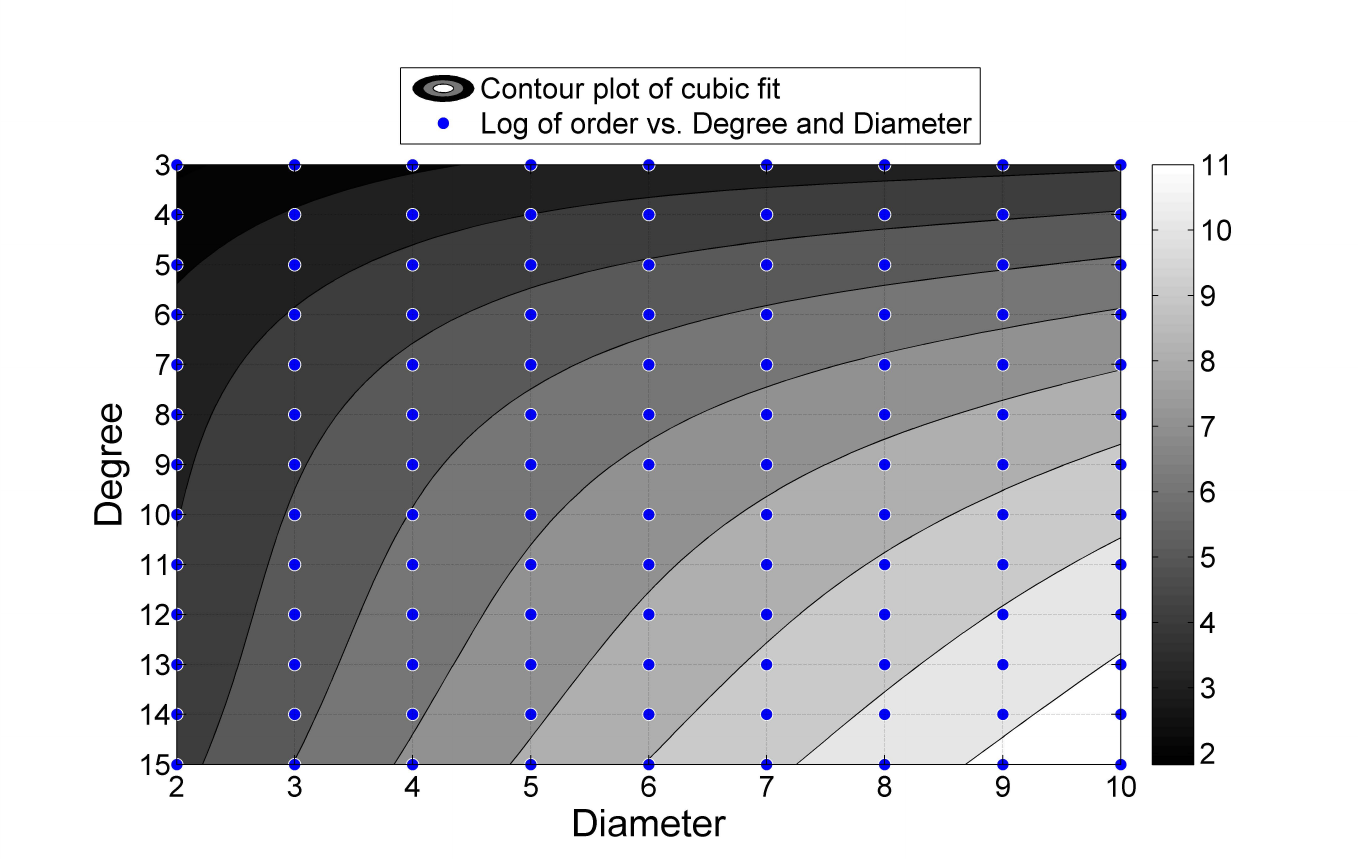}
\end{center}
\caption{Contour plot of the bi-cubic polynomial $f(\Delta, D)$, compared to the logarithmic Degree$/$Diameter table 
}
\label{fig:cubic-contour1B}
\end{figure}

Figure \ref{fig:cubic-contour1B} shows the approximation of the bi-cubic polynomial $f(\Delta, D)$ to the logarithms of the graph orders. 


\begin{figure}[htbp]
\begin{center}
	 	\includegraphics[width=0.9\textwidth]{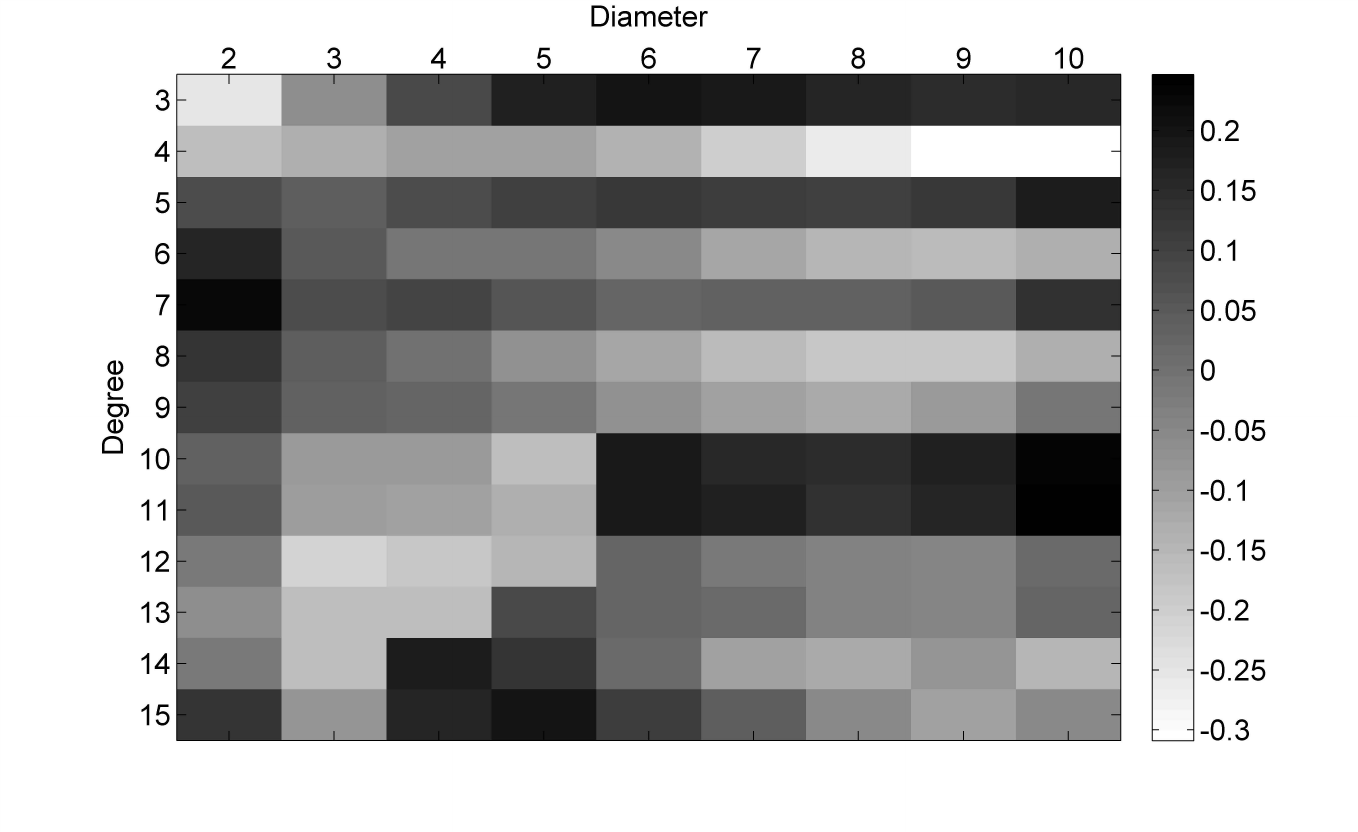}
\end{center}
\caption{Predicted (logarithmic) table values minus actual (logarithmic) values. Light entries correspond to negative differences. Darker colors correspond to positive differences.}
\label{fig:TableBW-cubic-15B}
\end{figure}

Fitting a function to the table provides insight in several ways. Obviously, it gives an idea on how the table grows, and how effective  the different construction methods are. Second, from this approximation we can figure out which are the table entries where an improvement is more likely to occur in the near future, because their current value lies beneath the value predicted by the function $e^{f(\Delta, D)}$. Figure \ref{fig:TableBW-cubic-15B} summarizes this analysis in a graphical way. Entry $(\Delta, D)$ contains a color representation of the predicted value $f(\Delta, D)$ minus the actual value of the logarithmic $(\Delta, D)$ table. The lighter the color, the better the actual value is, in comparison with the value predicted by $f(\Delta, D)$. Dark entries correspond to \lq bad\rq \ values, which lie below their prediction. 

There are several interesting things that pop out immediatlely: First of all, odd rows tend to be darker than even ones, especially for lower degrees. However, this should be no surprise to us, since the upper bounds show the same pattern, as we saw in Section \ref{sec:def}. When we look at the percentages (next), we see that even rows actually do not fare better than odd ones. 


The best opportunities for improvement seem to lie in rows $10$ and $11$, from diameters $6$ to $10$, and also in the small black region of rows $14, 15$, columns $4, 5$. On the other hand we get a light triangular area in the low diameters, in rows $10 - 15$. That is where we have been able to perform a more detailed search with Algorithm \ref{alg:search}. 

\begin{figure}[htbp]
\begin{center}
	 	\includegraphics[width=0.9\textwidth]{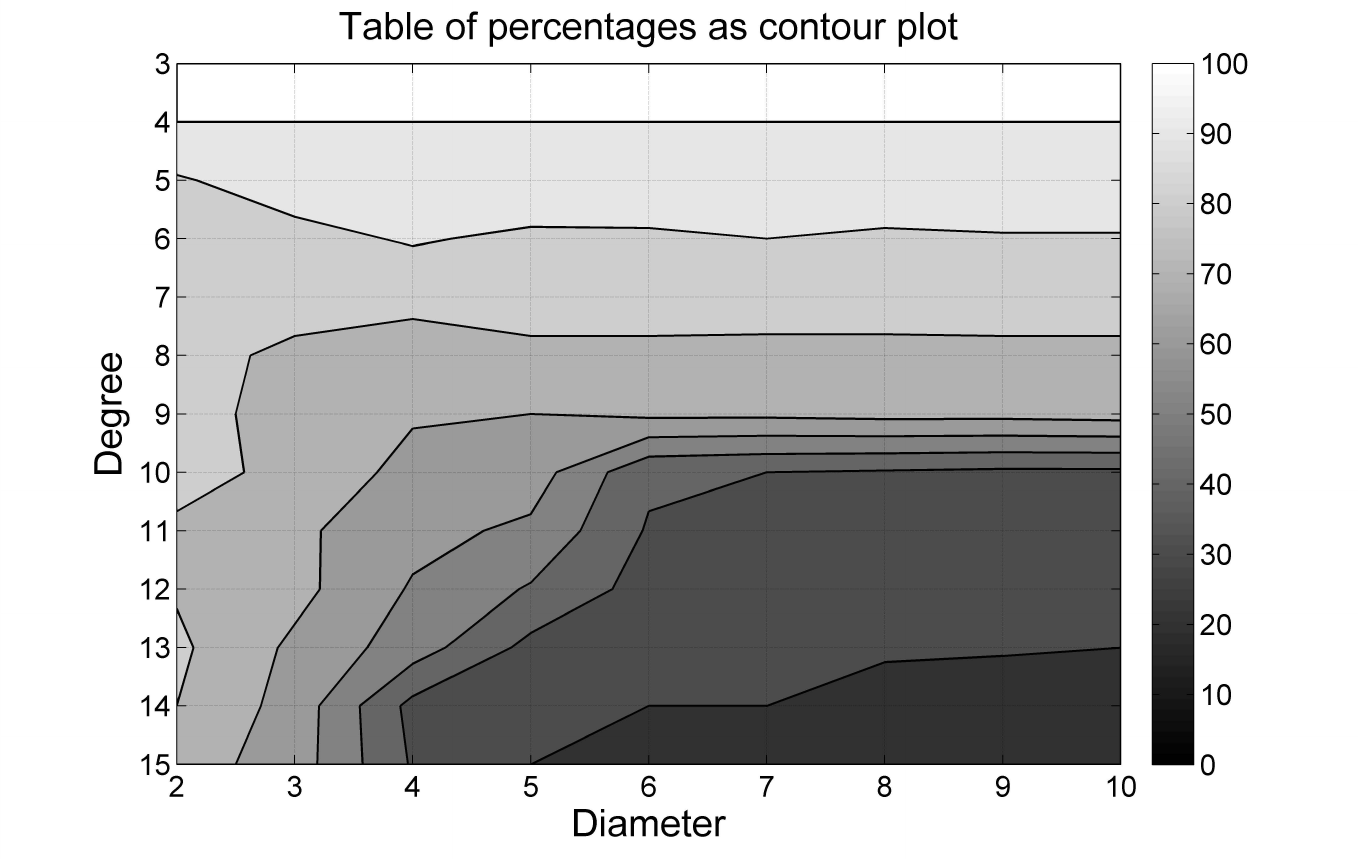}
\end{center}
\caption{Table of percentages shown as a contour plot}
\label{fig:ContourActualPercentages1}
\end{figure}

Another important analysis tool is given by the table of percentages, which indicates for each $\Delta$ and $D$ the percentage of the upper bound that has been achieved by the largest known graph. Figure \ref{fig:ContourActualPercentages1} shows the table of percentages as a contour plot. 

\begin{figure}[htbp]
\begin{center}
	 	\includegraphics[width=0.9\textwidth]{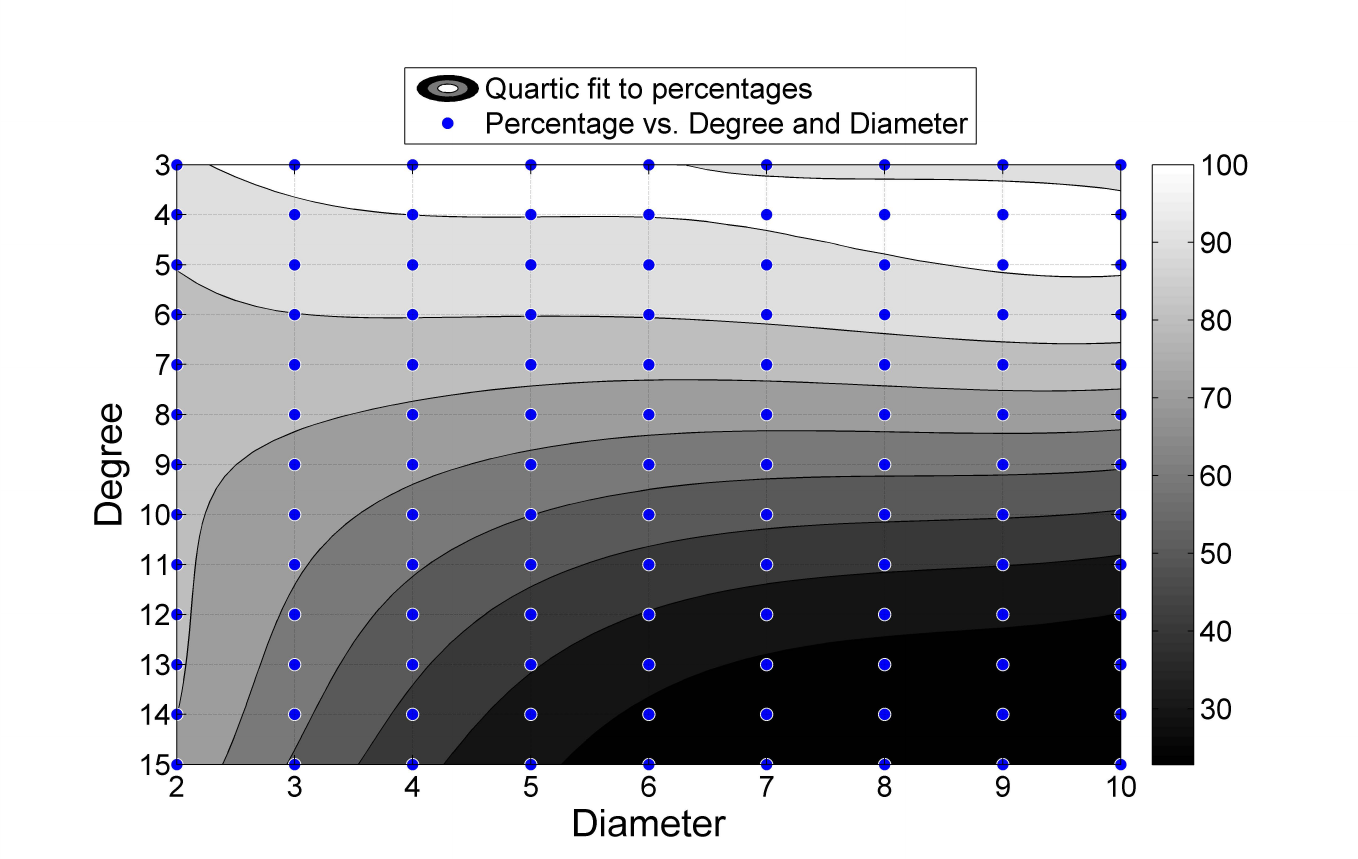}
\end{center}
\caption{Contour plot of $q(\Delta,D)$, which approximates the table of percentages}
\label{fig:Quartic-percentage-fit1A}
\end{figure}

It is not unreasonable to assume that the percentages should behave smoothly, therefore we can subject the table of percentages to the same analysis that we have used for the table of upper bounds and the table of graph orders. In this case, a bi-cubic polynomial is not flexible enough to approximate the table of percentages, therefore, we have fitted the bi-quartic polynomial 

\begin{equation}
\label{eq:fit3}
\begin{split}
q(\Delta,D) = & \; 98.77 -13.16 \Delta + 17.23 D + 1.214 \Delta^2 + 4.512 \Delta D \\ 
              &-6.777 D^2 + 0.03356 \Delta^3 -0.89 \Delta^2 D + 0.2302 \Delta D^2 \\
              &+ 0.6634 D^3 -0.003327 \Delta^4 + 0.02822 \Delta^3 D \\
              &+ 0.00937 \Delta^2 D^2 -0.01317 \Delta D^3 -0.02335 D^4
\end{split}
\end{equation}

to the table of percentages. The contour plot of $q(\Delta,D)$ is given in Figure \ref{fig:Quartic-percentage-fit1A}. 

\begin{figure}[htbp]
\begin{center}
	 	\includegraphics[width=0.9\textwidth]{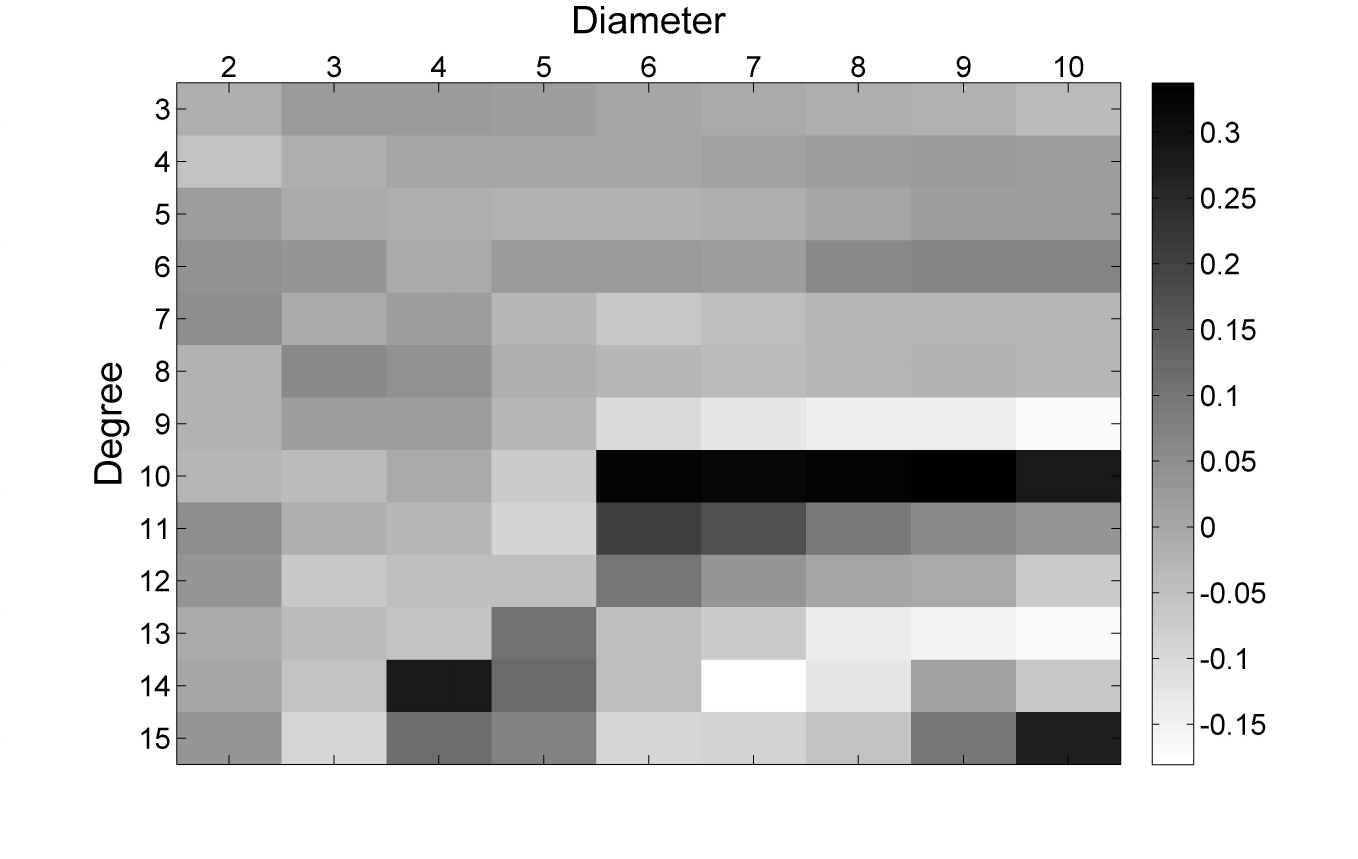}
\end{center}
\caption{Normalized differences between the predicted percentages and the actual percentages. Light entries correspond to negative differences. Darker colors correspond to positive differences.}
\label{fig:Percentage-differences-norm-BW1B}
\end{figure}

Now, for each $\Delta$ and $D$ we can construct the difference between the predicted value $q(\Delta,D)$ and the actual percentage. Figure \ref{fig:Percentage-differences-norm-BW1B} shows these differences, normalized by the actual percentages. Again, the lighter the color, the better the graph is, in comparison with its expected percentage. Figure \ref{fig:Percentage-differences-norm-BW1B} confirms the main conclusions drawn from the analysis of Figure \ref{fig:TableBW-cubic-15B}.

\section{Open problems}
\label{sec:open}

The apparent simplicity of circulant graphs is deceiving, and they are a reservoir of interesting open problems. The main theoretical question here has to do with finding sharper upper bounds for circulant graphs. The bounds $F(t,D)$ and $F'(t,D)$ given by Equations \ref{eq:circumoore}, \ref{eq:delannoy}, \ref{eq:oddbound}, and \ref{eq:evenbound} are clearly not sharp. All graphs in Table \ref{tab:circulant1}, from degrees 5 to 8 are known to be optimal, and yet they lie below the upper bounds. Some sharper lower bounds have been found for specific degrees or diameters, such as Eq. \ref{eq:triple}, but no general sharper lower bound is known up to date. 

The upper bounds $F(t,D)$ and $F'(t,D)$ have been obtained only by using the commutativity of the underlying group. However, circulant graphs are not even optimal in the class of abelian Cayley graphs (see \cite{Dou04} and also \cite{Fiol14}), despite some papers claiming so \cite{fallacy}. 

As for the lower bounds, there are few constructions for $\Delta > 8$. Our Cartesian product construction may clearly be generalized to obtain a general lower bound, which may then be compared analytically with other existing constructions. 

Regarding computer-based search, this problem poses several challenges. The search space grows very quickly with $\Delta$ and $D$, making a systematic search unfeasible. This problem could be alleviated by combining the search procedure with some isomorph rejection sieve, or some powerful heuristics. Additionally, other computer-based methods in the style of \cite{Mona02} could be devised, to obtain families of large circulant graphs with low diameter. 

In this paper we have only addressed the construction of undirected circulants, but in principle, the methods described here can be extended to directed, or even mixed graphs. 


\section*{Acknowledgements}
\label{sec:acknow}

We feel indebted to Charles Delorme for sharing his unpublished results with us. 



\end{document}